\documentclass[12pt,a4paper]{article}
\usepackage{amsmath,amsthm,amssymb,latexsym}
\usepackage{enumerate}
\usepackage[latin1]{inputenc}
\usepackage[OT2,OT1]{fontenc}
\usepackage{xy}
\xyoption{all}

\def\cyr{\fontencoding{OT2}\fontfamily{wncyr}\selectfont}

\newtheorem{theorem}{Theorem}
\newtheorem{lemma}[theorem]{Lemma}
\newtheorem{corollary}[theorem]{Corollary}
\newtheorem{proposition}[theorem]{Proposition}

\renewcommand{\phi}{\varphi}

\newcommand{\ran}{\operatorname{im}}
\newcommand{\im}{\operatorname{im}}
\newcommand{\rank}{\operatorname{rank}}
\newcommand{\str}{\operatorname{strk}}

\begin{document}

\title{Isolated subsemigroups in the variants of $\mathcal{T}_n$}
\author{Volodymyr Mazorchuk and Galyna Tsyaputa}
\date{}
\maketitle

\begin{abstract}
We classify all isolated, completely isolated and convex subsemigroups in 
the semigroup $\mathcal{T}_n$ of all transformations of an $n$-element set,
considered as the semigroup with respect to a sandwich operation.
\end{abstract}

\section{Introduction and description of the results}\label{s1}

For a semigroup, $(S,\cdot)$, and an element, $a\in S$, one can consider the 
{\em variant} $(S,*_a)$ of $S$, for which the {\em deformed} or {\em sandwich} 
multiplication $*_a:S\times S\to S$ with the sandwich element $a$ is defined 
as follows:  $x*_ay=xay$ for all $x,y\in S$. This construction was proposed in 
\cite{Ly} and studied later on for various classes of semigroups by several 
authors, see for example \cite{Ch1,Ch2,Hi1,Hi2,Hu,MMT,Su,Sy,Ts1,Ts2,Ts3,KT} and 
others. One of the main motivations for this study is the fact that usually the 
variants of the semigroups have more interesting and richer structure than the 
original semigroups.

The basic question, which one faces in the study of  variants
of the semigroup $S$ is: {\em For which $a,b\in S$ are the variants
$(S,*_a)$ and $(S,*_b)$ isomorphic?} For such classical transformation
semigroups on a given set, $N$, as the full symmetric inverse semigroup 
$\mathcal{IS}(N)$, the full transformation semigroup $\mathcal{T}(N)$, and
the semigroup $\mathcal{PT}(N)$ of all partial transformations, the 
isomorphism criteria for their variants were obtained in
\cite{Sy,Ts1,Ts3}. In particular, in all the above cases it was shown that 
(up to isomorphism) one can always assume that $a$ is an idempotent of $S$.

A natural class of subsemigroups of the given semigroup $S$ is formed by
the so-called {\em isolated} subsemigroups. A subsemigroup, $T$, of $S$ 
is called {\em isolated} provided that the following condition is satisfied: 
$x^n\in T$ for some $n\in\mathbb{N}$ implies $x\in T$ for all $x\in S$, 
see for example \cite{M-S}. This class
contains various subclasses, for example the class of completely isolated 
subsemigroups and some others, see Section~\ref{s2} or \cite{M-S} for details.
Isolated subsemigroups of some classical semigroups are known. 
For example,  a complete classification  of isolated and 
completely isolated subsemigroups of $\mathcal{IS}(N)$ in the case of finite
$N$ was obtained in \cite{GK}. For the factor power $\mathcal{FP}(\mathcal{S}_n)$
of the symmetric group $\mathcal{S}_n$ the same problem was solved in \cite{GM}, 
and for the Brauer semigroup $\mathfrak{B}_n$ in \cite{Ma}. 

In the present paper we present a complete description of the isolated 
subsemigroups for all (up to isomorphism) variants of the full transformation 
semigroup $\mathcal{T}(N)$ in the case of a finite set $N$. We also classify
the classes of completely isolated, right convex, left convex and convex 
subsemigroups, see Section~\ref{s2} for the corresponding definitions. 

Let us now briefly describe 
our approach, results and the structure of the paper. The results of
\cite{GK,GM,Ma} suggest to split the study of the isolated subsemigroups into
three different cases. In \cite{GK,GM,Ma} these are: the case of 
invertible elements, the case of maximal subgroups corresponding to  
maximal non-invertible idempotents, and the rest. For the variants of
semigroups the above cases should first be adjusted as all non-trivial 
variants never contain any unit and, as a consequence, any invertible elements. 
However, in the case of transformation semigroups there is a natural 
``substitution'' for this. For  transformation semigroups there is a natural
notion of the rank (the cardinality of the image) and
stable rank (the rank of the unique
idempotent in the cyclic subsemigroup, generated by the given element).
Using these notions we split our study into the following three cases:
the case of the maximal possible stable rank $l$ (which equals the rank
of the sandwich element), the case of the stable rank $l-1$, and
the rest. Even in the first case we obtain quite a striking difference
with the classical case. In the classical case all invertible elements obviously
form a minimal isolated subsemigroup. In contrast to this, 
for every variant of $\mathcal{T}_n$ we
show in Subsection~\ref{s4.1} that the set $\mathcal{T}_n^{(l)}$ of all 
elements of stable rank $l$ is a subsemigroup, and that there is a natural 
congruence on $\mathcal{T}_n^{(l)}$, the quotient modulo which is a rectangular 
band. Moreover, the natural projection onto the quotient induces a bijection 
between the set of all (isolated) subsemigroups of this rectangular band (which is
relatively easy to describe) and the set of all isolated subsemigroups of
$\mathcal{T}_n^{(l)}$. 

Not completely unexpected, but it turned out that the most interesting and the
most complicated case of our study is the case of the elements of stable 
rank $l-1$ (in the classical case certain subsemigroups, consisting of elements
of rank $l-1$, make the difference between the classes of isolated and completely
isolated subsemigroups). Such subsemigroups are studied
in Subsection~\ref{s4.2}, which is the heart of our paper. We show that 
the isolated subsemigroups of this form split into three classes. In every
class the semigroups are indexed by several discrete parameters which are
either elements or subsets of some (not very complicated) sets. The structure
of the simplest of these three classes is quite similar to the structure of
isolated subsemigroups  in the case of $\mathcal{T}_n^{(l)}$. 
Two other classes are in some sense ``composed'' from 
several copies of the first class, which leads to quite nice but non-trivial
description of isolated subsemigroups in the case of stable rank $l-1$.

The last case, that is the case of stable rank at most $l-2$, 
is studied in Subsection~\ref{s4.3}. This case is relatively ``poor'' with
respect to the structure of the isolated subsemigroups. In fact we show that
if an isolated subsemigroup contains an element, whose stable rank does 
not exceed $l-2$, then it contains all elements, whose  stable ranks do not 
exceed $l-1$ (and the last $l-1$ is not a typo, it is indeed $l-1$ and
not $l-2$). In Subsection~\ref{s4.4}
we present our main theorem about the classification of all isolated
subsemigroups in all variants of $\mathcal{T}_n$. In Section~\ref{s5}
we use this theorem to classify all completely isolated and
convex subsemigroups. We finish the paper with deriving the corollaries 
from our results for the exceptional cases (that is the cases when the sandwich 
element has the maximal possible rank $n$, the minimal possible rank $1$, or 
rank $2$) in Section~\ref{s6}.

The basic generalities on isolated subsemigroups are collected in 
Section~\ref{s2}. In Section~\ref{s3} we fix the notation and collect
some generalities on the semigroup $\mathcal{T}_n$ and its variants.
We would like to emphasize once more that  we always assume the sandwich 
element to be an idempotent, which we can do up to isomorphism of the sandwich 
semigroup $\mathcal{T}_n$. Among other basic facts, in Subsection~\ref{s3.3} we 
describe all idempotents and the corresponding maximal subgroups in all variants of 
$\mathcal{T}_n$. 

\section{Generalities on isolated subsemigroups}\label{s2}

In this section we collect some generalities on isolated and completely isolated
subsemigroups, which we failed to find in the literature and which we think might
be of independent interest. Some of these basic properties will be often used
in the paper without reference.

Let $S$ be a semigroup. A subsemigroup, $T\subset S$, is called
\begin{itemize}
\item {\em isolated} provided that for all $x\in S$ the condition
$x^n\in T$ for some $n>0$ implies $x\in T$;
\item {\em completely isolated} provided that $xy\in T$ implies $x\in T$ or
$y\in T$ for all $x,y\in S$;
\item {\em right convex} provided that $xy\in T$ implies $y\in T$ for all 
$x,y\in S$;
\item {\em left convex} provided that $xy\in T$ implies $x\in T$ for all 
$x,y\in S$;
\item {\em convex} provided that $xy\in T$ implies $x\in T$ and
$y\in T$ for all $x,y\in S$.
\end{itemize}
We refer the reader to \cite{CP,M-S} for details on these definitions.
Denote by $\mathfrak{I}(S)$, $\mathfrak{CI}(S)$, $\mathfrak{RC}(S)$, 
$\mathfrak{LC}(S)$ and $\mathfrak{C}(S)$, the sets of all isolated,
completely isolated, right convex, left convex or convex subsemigroups of
$S$ respectively. For $X\subset S$ we denote by $\overline{X}$ the complement of
$X$ in $S$.

From the definition one immediately obtains the following properties
(here everywhere $T$ is assumed to be a subsemigroup of $S$): 
\begin{enumerate}[(I)]
\item\label{br1} $\mathfrak{I}(S)$ is closed with respect to 
(non-empty) intersections, in particular, for every $x\in S$ there exists the 
unique isolated subsemigroup $\mathtt{I}(x)$ of $S$ such that 
$x\in \mathtt{I}(x)$ and which is minimal with respect to this property.
\item\label{br2} $T\in \mathfrak{I}(S)$ if and only if $\overline{T}$ is either empty
or a (not necessarily disjoint) union of subsemigroups of $S$.
\item\label{br3} $\mathfrak{CI}(S)$ is closed with respect to taking the (non-empty) 
complement in $S$.
\item\label{br4} $T\in \mathfrak{CI}(S)$ if and only if $\overline{T}$ is either empty or
a subsemigroup of $S$.
\item\label{br5} $T\in \mathfrak{RC}(S)$ if and only if $\overline{T}$ is either empty or
a left ideal.
\item\label{br6} $T\in \mathfrak{LC}(S)$ if and only if $\overline{T}$ is either empty or
a right ideal.
\item\label{br7} $T\in \mathfrak{C}(S)$ if and only if $\overline{T}$ is either empty or a (two-sided) ideal.
\item\label{br8} If $S_1$ and $S_2$ are two semigroups and $f:S_1\to S_2$ is a 
homomorphism then taking the pre-image (provided that it is non-empty) preserves the 
properties of being isolated, completely isolated, left convex, right convex and convex, 
that is, it induces a map from $\mathfrak{I}(S_2)$ to $\mathfrak{I}(S_1)$ and so on.  
\item
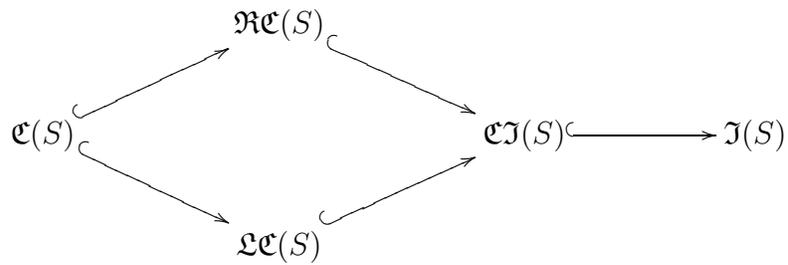
\label{br9} 
\begin{displaymath}
\xymatrix{
 && \mathfrak{RC}(S)\ar@{^{(}->}[rrd]&& &&  \\
\mathfrak{C}(S)\ar@{^{(}->}[rru]\ar@{^{(}->}[rrd] && && 
\mathfrak{CI}(S)\ar@{^{(}->}[rr]&& \mathfrak{I}(S) \\
 && \mathfrak{LC}(S)\ar@{^{(}->}[rru]&& &&  
}
\end{displaymath}
\end{enumerate}

We remark that for $T_1,T_2\in \mathfrak{CI}(S)$ the subsemigroup $T_1\cap T_2$
does not belong to $\mathfrak{CI}(S)$ in general. An example can be found in
Section~\ref{s5}. From Section~\ref{s5} it also follows that all the five
classes above can be different. Now we would like to characterize some of the 
above subsemigroups in terms of certain homomorphisms and congruences on $S$.

\begin{proposition}\label{pr1}
\begin{enumerate}[(a)]
\item\label{pr1.1} $T\in \mathfrak{C}(S)$ if and only if there exists 
$f:S\to(\mathbb{Z}_2,\cdot)$ such that $T=f^{-1}(1)$.
\item\label{pr1.2} Let $\mathfrak{C}(S)=\{S,T_1,\dots,T_k\}$. Then 
$\{X_1\cap\dots\cap X_k\neq \varnothing: X_i=T_i\text{ or }X_i=\overline{T}_i\}$
defines a congruence, $\rho$, on $S$, the quotient modulo which embeds into
$(\mathbb{Z}_2,\cdot)^k$. 
\item\label{pr1.3} For $x\in S$ the congruence class $\rho_x$ of $\rho$,
containing $x$, is an isolated subsemigroup of $S$, in particular, the congruence
$\rho$ gives rise to a decomposition of $S$ into a disjoint union of isolated
subsemigroups.
\end{enumerate}
\end{proposition}

\begin{proof}
Let $T$ be a proper convex subsemigroup of $S$. Consider the equivalence
relation on $S$ with two equivalence classes: $T$ and $\overline{T}$. This
relation is a congruence on $S$ and the corresponding quotient is 
isomorphic to $(\mathbb{Z}_2,\cdot)$ (the isomorphism is given by
$T\mapsto 1$, $\overline{T}\mapsto 0$).
Conversely, note that $\{1\}\in \mathfrak{C}((\mathbb{Z}_2,\cdot))$.
Hence if $f:S\to(\mathbb{Z}_2,\cdot)$ then \eqref{br8} implies that 
$f^{-1}(1)\in \mathfrak{C}(S)$, which proves \eqref{pr1.1}.

For $i=1,\dots,k$ let $f_i:S\rightarrow \mathbb{Z}_2$ denote the
homomorphism such that $f_i(T_i)=1$ and $f_i(\overline{T}_i)=0$. 
Then by \eqref{pr1.1} the kernel of the homomorphism $(f_1,\dots,f_k)$ 
is just $\rho$. This proves \eqref{pr1.2}.

\eqref{pr1.3} follows from \eqref{br1}, \eqref{br9} and the fact that
$T_i,\overline{T}_i\in \mathfrak{CI}(S)$ for all $i$.  
\end{proof}

\begin{proposition}\label{pr3}
Let $\mathfrak{K}$ be one of the classes $\mathfrak{I}(S)$,
$\mathfrak{CI}(S)$, $\mathfrak{RC}(S)$, $\mathfrak{LC}(S)$ or
$\mathfrak{C}(S)$. If $T$ is a subsemigroup of $S$ and is a 
(not necessarily disjoint) union of subsemigroups from $\mathfrak{K}$,
then $T\in \mathfrak{K}$.
\end{proposition}

\begin{proof}
Follows from the fact that the conditions on $\overline{T}$, which
appear in \eqref{br2}, \eqref{br4}, \eqref{br5}, \eqref{br6}, and \eqref{br7},
are closed with respect to  intersections.
\end{proof}

We would like to remark that the union of isolated subsemigroups 
(or any other type of subsemigroups defined above) does not need to be a 
semigroup in general. One easily constructs examples using our results, 
presented in Section~\ref{s4} and Section~\ref{s5}.

Denote by $E(S)$ the set of all idempotents of $S$ and for $e\in E(S)$
define 
\begin{displaymath}
\sqrt{e}=\{x\in S: x^m=e\text{ for some } m>0\}.
\end{displaymath}

\begin{proposition}\label{pr4}
If $\sqrt{e}$ is a subsemigroup of  $S$ then 
$\sqrt{e}$ is a minimal (with respect to inclusions) isolated
subsemigroup of $S$, in particular, $\sqrt{e}=\mathtt{I}(x)$ 
for every $x\in \sqrt{e}$.
\end{proposition}

\begin{proof}
Assume that $\sqrt{e}$ is a subsemigroup of  $S$. If
$y^m\in \sqrt{e}$ for some $m>0$ then $(y^m)^k=e$ for some $k>0$ and hence 
$y\in \sqrt{e}$ implying that $\sqrt{e}$ is isolated.

Let now $x\in \sqrt{e}$. Then $\mathtt{I}(x)$ contains all powers of
$x$, in particular, it contains $e$. Since $\mathtt{I}(x)$ is isolated, it
follows that $\mathtt{I}(x)$ must contain $\sqrt{e}$. However, as we
have already shown, $\sqrt{e}$ is isolated itself. Hence
$\sqrt{e}=\mathtt{I}(x)$ by the minimality of $\mathtt{I}(x)$. This 
completes the proof.
\end{proof}

\section{$\mathcal{T}_n$ and its variants}\label{s3}

\subsection{Notation}\label{s3.1}

Throughout the paper  we fix a positive integer, $n$, and set
$N=\{1,2,\dots,n\}$. We denote by $\mathcal{T}_n$ the semigroup of
all maps $\beta:N\to N$ with respect to the composition $\cdot$ of 
maps from the left to the right (in contrast with the standard 
composition $\circ$ which works from the right to the left). We 
adopt the standard notation $\beta(x)$ for the value of
$\beta$ on the element $x\in N$.
In particular, we always have $\beta\gamma(x)=\gamma(\beta(x))$
for all $\beta,\gamma\in \mathcal{T}_n$ and $x\in N$.

For $\beta\in\mathcal{T}_n$ we denote by $\ran(\beta)$ the {\em image}
of $\beta$, that is the set $\{\beta(x):x\in N\}$, and define the
{\em rank} of $\beta$ via $\rank(\beta)=|\ran(\beta)|$. Define also
the equivalence relation $\Lambda_{\beta}$ on $N$ as follows:
$x\Lambda_{\beta} y$, $x,y\in N$,  if and only if $\beta(x)=\beta(y)$.
Sometimes it will be convenient to consider $\Lambda_{\beta}$
as an unordered (disjoint) union of equivalence classes and use the notation
\begin{displaymath}
\Lambda_{\beta}=B_1^{\beta}\cup B_2^{\beta}\cup\dots \cup B_k^{\beta}
\end{displaymath}
if  $\rank(\beta)=k$. If the element $\beta$ is clear from the context, we
sometimes may skip it as the upper index.  The relation $\Lambda_{\beta}$
is called the {\em kernel} of $\beta$.

If $\beta\in \mathcal{T}_n$ and $N=N_1\cup N_2\cup\dots\cup N_k$ is a 
decomposition of $N$ into a disjoint union of subsets such that 
for all $i=1,\dots,k$ we have $\beta(x)=\beta(y)=b_i$ for all $x,y\in N_i$
and for some $b_i\in N$, we will use the following notation for $\beta$:
\begin{displaymath}
\beta=\left(
\begin{array}{cccc}
N_1 & N_2 & \dots & N_k\\
b_1 & b_2 & \dots & b_k
\end{array}
\right).
\end{displaymath}
Note that we do not require $\rank(\beta)=k$, which means that some of 
$b_i$'s may coincide in general.

For $i\in N$ we denote by $\theta_i$ the element of $\mathcal{T}_n$
defined via $\theta_i(x)=i$ for all $x\in N$. Note that $\theta_i$ is a
right zero for $\mathcal{T}_n$. For a set, $X$, we denote by 
$\mathcal{S}(X)$ the symmetric group on $X$.

\subsection{Variants of $\mathcal{T}_n$ and their isomorphism}\label{s3.2}

For $\alpha\in \mathcal{T}_n$ let $(\mathcal{T}_n,*_{\alpha})$ be the
{\em variant} of $\mathcal{T}_n$, that is the semigroup $\mathcal{T}_n$ with
the {\em sandwich operation} $\beta*_{\alpha}\gamma=\beta\alpha\gamma$, 
$\beta,\gamma\in \mathcal{T}_n$, with respect to the {\em sandwich element} 
$\alpha$. The following statement is proved in \cite{Sy,Ts1}:

\begin{theorem}\label{t5}
Let $\alpha_1,\alpha_2\in \mathcal{T}_n$. Then the following statements are
equivalent:
\begin{enumerate}
\item\label{t3.1} $(\mathcal{T}_n,*_{\alpha_1})\cong (\mathcal{T}_n,*_{\alpha_2})$.
\item\label{t3.2} $\alpha_1\sigma=\tau\alpha_2$ for some  $\sigma,\tau\in\mathcal{S(N)}$.
\item\label{t3.3} 
For every $i=1,\dots,n$ the decompositions $\Lambda_{\alpha_1}$ and $\Lambda_{\alpha_2}$ 
contain the same number of blocks of cardinality $i$. 
\end{enumerate}
\end{theorem}

From Theorem~\ref{t5} it follows that every variant of $\mathcal{T}_n$ is isomorphic
to a variant, the sandwich element of which is an idempotent. Hence from now on we
fix $\alpha\in E(\mathcal{T}_n)$ of rank $l$. Set $*=*_{\alpha}$ and let
\begin{displaymath}
\alpha=\left(
\begin{array}{cccc}
A_1 & A_2 & \dots & A_l\\
a_1 & a_2 & \dots & a_l
\end{array}
\right).
\end{displaymath}
Then $\Lambda_{\alpha}=A_1\cup A_2\cup\dots\cup A_l$, $\ran(\alpha)=\{a_1,\dots,a_l\}$,
and $a_i\in A_i$ for all $i$. We also set $A=\{a_1,\dots,a_l\}$.

Let $\beta\in\mathcal{T}_n$ and $\varepsilon\in E(\mathcal{T}_n,*)$ be the unique 
idempotent in the cyclic subsemigroup generated by $\beta$. We define the
{\em stable rank} of $\beta$ as $\str(\beta)=\rank(\varepsilon)$, and set 
$\Lambda_{\beta}^{\mathrm{st}}=\Lambda_{\varepsilon}$,
$\ran^{\mathrm{st}}(\beta)=\ran(\varepsilon)$. Note that
$\Lambda_{\beta}^{\mathrm{st}}=\Lambda_{\beta\alpha}^{\mathrm{st}}$ and
$\ran^{\mathrm{st}}(\beta)=\ran^{\mathrm{st}}(\alpha\beta)$ for all
$\beta\in (\mathcal{T}_n,*)$. For $\beta\in \mathcal{T}_n$
and $m\in\mathbb{N}$ we set 
\begin{displaymath}
\beta^m=\underbrace{\beta\beta\dots\beta}_{m\,\,\mathrm{times}}\quad\quad
\text{and}\quad\quad
\beta^{*m}=\underbrace{\beta*\beta*\dots*\beta}_{m\,\,\mathrm{times}}.
\end{displaymath}
For $i=1,\dots,l$ we denote by $\mathcal{T}_n^{(i)}$ the set of all elements in 
$(\mathcal{T}_n,*)$ of stable rank $i$. 

For some other properties of $(\mathcal{T}_n,*)$, for example for the description
of Green's relations or the automorphism group, we refer the reader to \cite{Sy,Ts2}.

\subsection{Description of idempotents and maximal subgroups}\label{s3.3}

To proceed with the study of isolated subsemigroups of $(\mathcal{T}_n,*)$
we need to describe the idempotents in this semigroup.

\begin{theorem}\label{t6}
\begin{enumerate}[(a)]
\item\label{t6.1} Let 
\begin{equation}\label{idemp}
\varepsilon=\left(
\begin{array}{cccc}
E_1 & E_2 & \dots & E_k\\
e_1 & e_2 & \dots & e_k
\end{array}
\right)\in \mathcal{T}_n
\end{equation}
be an element of rank $k$. Then
$\varepsilon\in E(\mathcal{T}_n,*)$ if and only if there exists an injection,
$f:\{1,\dots,k\}\to \{1,\dots,l\}$, such that the following two conditions are 
satisfied:
\begin{enumerate}[(i)]
\item\label{t6.1.1} $e_i\in A_{f(i)}$ for all $i=1,\dots,k$,
\item\label{t6.1.2} $a_{f(i)}\in E_i$ for all $i=1,\dots,k$.
\end{enumerate}
\item\label{t6.2} 
\begin{displaymath}
|E(\mathcal{T}_n,*)|=
\sum_{\varnothing\neq X\subset\{1,\dots,l\}}
\left(\prod_{i\in X}|A_i|\right)\cdot |X|^{n-|X|}.
\end{displaymath}
\end{enumerate}
\end{theorem}

\begin{proof}
If $\varepsilon\in E(\mathcal{T}_n,*)$ then
$\varepsilon=\varepsilon*\varepsilon=\varepsilon\alpha\varepsilon$. Therefore
$\rank(\varepsilon)=\rank(\varepsilon\alpha)$ and hence all $e_i$'s belong to
different blocks of $\Lambda_{\alpha}$. Define $f:\{1,\dots,k\}\to \{1,\dots,l\}$
via the requirement $e_i\in A_{f(i)}$ for all $i=1,\dots,k$ and \eqref{t6.1.1} follows.
On the other hand, let $\varepsilon(a_{f(i)})=e_j$. Then
\begin{displaymath}
\varepsilon(a_{f(i)})=\varepsilon\alpha\varepsilon(a_{f(i)})=
\alpha\varepsilon(e_j)=\varepsilon(a_{f(j)}),
\end{displaymath}
which implies $f(i)=f(j)$ and thus $a_{f(i)}\in E_i$ for all $i=1,\dots,k$
giving \eqref{t6.1.2}. Conversely, if both \eqref{t6.1.1} and \eqref{t6.1.2} are
satisfied, a direct calculation shows that $\varepsilon\in E(\mathcal{T}_n,*)$.
This proves \eqref{t6.1}.

To prove \eqref{t6.2} we count the number of idempotents, for which
$\ran(f)=X\subset\{1,\dots,l\}$, where $f$ is taken from \eqref{t6.1}. 
If $X$ is fixed, we should independently choose elements $e_i\in A_i$ for all
$i\in X$, and then every element from the set $N\setminus\{a_{f(i)}:i\in X\}$ 
should be associated to some block of $\Lambda_{\varepsilon}$.
The formula of \eqref{t6.2} now follows by an application of the product rule and
the sum rules.
\end{proof}

The above description of idempotents allows us to describe the maximal subgroups
in $(\mathcal{T}_n,*)$:

\begin{corollary}\label{c7}
Let $\varepsilon$ be an idempotent in $(\mathcal{T}_n,*)$ of rank $k$, defined
via \eqref{idemp}. Then the corresponding maximal subgroup $G(\varepsilon)$ of 
$(\mathcal{T}_n,*)$ is isomorphic to the symmetric group $\mathcal{S}(\{e_1,\dots,e_k\})$
and consists of all $\beta\in \mathcal{T}_n$ which satisfy the following two conditions:
\begin{enumerate}[(i)]
\item\label{c7.1} $\Lambda_{\beta}=\Lambda_{\varepsilon}$,
\item\label{c7.2} $\ran(\beta)=\ran(\varepsilon)$.
\end{enumerate}
\end{corollary}

\begin{proof}
$\beta\in G(\varepsilon)$ if and only if $\beta\in\sqrt{\varepsilon}$ and
$\varepsilon*\beta=\beta*\varepsilon=\beta$. In particular, it follows
that $\rank(\beta)=\rank(\varepsilon)$. For such $\beta$ the condition
\eqref{c7.1} is equivalent to $\varepsilon*\beta=\beta$ and \eqref{c7.2} is
equivalent to $\beta*\varepsilon=\beta$. This completes the proof.
\end{proof}

\subsection{Some homomorphisms from $(\mathcal{T}_n,*)$ to $\mathcal{T}_n$
and $\mathcal{T}(A)$}\label{s3.4}

In this subsection we construct two homomorphisms from 
$(\mathcal{T}_n,*)$ to the usual $\mathcal{T}_n$ and an epimorphism
to $\mathcal{T}_m$ for some $m\leq n$. For every $\beta\in\mathcal{T}_n$ set
\begin{displaymath}
\varphi_{l}(\beta)=\alpha\beta,\quad\quad
\varphi_{r}(\beta)=\beta\alpha.
\end{displaymath}

\begin{lemma}\label{l15}
\begin{enumerate}[(a)]
\item\label{l15.1}
Both $\varphi_{l}$ and $\varphi_{r}$ are homomorphisms from
$(\mathcal{T}_n,*)$ to $\mathcal{T}_n$.
\item\label{l15.2} 
For every $\beta\in \mathcal{T}_n$ we have
\begin{displaymath}
\str(\beta)=\str(\varphi_{l}(\beta))=
\str(\varphi_{r}(\beta)).
\end{displaymath}
\end{enumerate}
\end{lemma}

\begin{proof}
We prove both statement for, say, $\varphi_{r}$. In other cases the arguments
are similar. For $\beta,\gamma\in \mathcal{T}_n$
we have:
\begin{displaymath}
\varphi_{r}(\beta*\gamma)=(\beta*\gamma)\alpha=
\beta\alpha\gamma\alpha=(\beta\alpha)(\gamma\alpha)=
\varphi_{r}(\beta)\varphi_{r}(\gamma),
\end{displaymath}
which proves \eqref{l15.1}.

We prove \eqref{l15.2} for $\varphi_{r}$ using the following argument.
First we observe  that 
$\beta*\beta=\beta\alpha*\beta$ as $\alpha\alpha=\alpha$.
Therefore, since
$\beta^{*(m+1)}=(\beta\alpha)^{*m}\beta$ for all $m$, we have
$\rank(\beta^{*(m+1)})\leq \rank((\beta\alpha)^{*m})$.
As $(\beta\alpha)^{*(m+1)}=\beta^{*(m+1)}\alpha$ for all $m$, we have
$\rank((\beta\alpha)^{*(m+1)})\leq \rank(\beta^{*(m+1)})$.
However, by the definition of the stable rank,
$\str(\beta)=\rank(\beta^{*m})$ and
$\str(\beta\alpha)=\rank((\beta\alpha)^{*m})$ for all $m$ big enough. 
This implies that $\str(\beta)=\str(\varphi_{r}(\beta))$.
For $\varphi_{l}$ the argument is analogous, which completes the proof.
\end{proof}

Let $\mathcal{T}(A)$ denote the full transformation semigroup on $A$ and observe
that $\alpha\beta(A)\subset A$ for all $\beta\in\mathcal{T}_n$. Hence
from Lemma~\ref{l15}\eqref{l15.1} we immediately obtain:

\begin{corollary}\label{chom}
The map $\overline{\varphi}:(\mathcal{T}_n,*)\to\mathcal{T}(A)$, defined via
\begin{displaymath}
\overline{\varphi}(\beta)=\varphi_{r}(\beta)|_{A},
\end{displaymath}
is an epimorphism.
\end{corollary}

Note that the congruence, which corresponds to $\overline{\varphi}$, 
was also used in \cite{Sy}.

\section{Classification of the isolated subsemigroups}\label{s4}

From the definition of $*$ it follows that the stable rank of an element
from $(\mathcal{T}_n,*)$ does not exceed $l$. The results of
\cite{GK,GM,Ma} suggest that for the description of the minimal isolated 
subsemigroups in $(\mathcal{T}_n,*)$ it would be convenient to consider 
three different cases, namely the elements of stable rank $l$, $l-1$, and at 
most $l-2$, separately.

\subsection{The case of stable rank $l$}\label{s4.1}

\begin{lemma}\label{l11}
\begin{enumerate}[(a)]
\item\label{l11.01} Let $\beta\in \mathcal{T}_n$. Then $\beta$ has stable
rank $l$ if and only if $\overline{\varphi}(\beta)$ is a bijection.
\item\label{l11.1} Let $\beta\in \mathcal{T}_n$ be an element of stable rank $l$.
\begin{displaymath}
\Lambda^{\mathrm{st}}_{\beta}=\Lambda_{\beta\alpha},\quad\quad
\text{ and }\quad\quad
\ran^{\mathrm{st}}(\beta)=\ran(\alpha\beta)=\beta(\{a_1,\dots,a_l\}).
\end{displaymath}
\item\label{l11.2} 
Let $\beta,\gamma\in \mathcal{T}_n$ be two elements of stable rank $l$. Then
\begin{displaymath}
\Lambda^{\mathrm{st}}_{\beta*\gamma}=\Lambda^{\mathrm{st}}_{\beta},\quad\quad
\text{ and }\quad\quad
\ran^{\mathrm{st}}(\beta*\gamma)=\ran^{\mathrm{st}}(\gamma).
\end{displaymath}
\end{enumerate}
\end{lemma}

\begin{proof}
If $\overline{\varphi}(\beta)$ is a bijection, then 
$\overline{\varphi}(\beta)^m$ is a bijection for any $m$. Hence the
stable rank of $\beta$ is at least $l$, which implies that it is exactly
$l$. On the other hand, the stable rank of $\beta$ does not exceed the
rank of $\overline{\varphi}(\beta)$, which means that the stable rank of 
$\beta$ is strictly smaller than $l$ in the case when 
$\overline{\varphi}(\beta)$ is not a bijection. This proves \eqref{l11.01}.

Since $\rank(\alpha)=l$ we have both that $|\ran(\alpha\beta)|\leq l$ and that
the number of blocks in $\Lambda_{\beta\alpha}$ does not exceed $l$. However,
since $\str(\beta)=l$ by our assumption, it follows that 
$|\ran(\alpha\beta)|=l$ and that the number of blocks in 
$\Lambda_{\beta\alpha}$ equals $l$. Since, obviously, 
$\Lambda_{\beta\alpha}\subset \Lambda^{\mathrm{st}}_{\beta}$ and
$\ran^{\mathrm{st}}(\beta)\subset \ran(\alpha\beta)$, we derive
$\Lambda^{\mathrm{st}}_{\beta}=\Lambda_{\beta\alpha}$ and
$\ran^{\mathrm{st}}(\beta)= \ran(\alpha\beta)$. The statement \eqref{l11.1} follows.

Since, obviously, $\Lambda_{\beta\alpha}\subset \Lambda^{\mathrm{st}}_{\beta*\gamma}$
and $\ran^{\mathrm{st}}(\beta*\gamma)\subset \ran(\alpha\gamma)$, 
the statement \eqref{l11.1} implies 
$\Lambda^{\mathrm{st}}_{\beta}\subset \Lambda^{\mathrm{st}}_{\beta*\gamma}$
and $\ran^{\mathrm{st}}(\beta*\gamma)\subset \ran^{\mathrm{st}}(\gamma)$. On the other  hand, as $\str(\beta)=\str(\gamma)=l$ the statement \eqref{l11.01}
implies that $\overline{\varphi}(\beta*\gamma)$ is a bijection and hence
$\str(\beta*\gamma)=l$ again by \eqref{l11.01}. 
The last equality, however, says that different blocks of
$\Lambda^{\mathrm{st}}_{\beta}$ should remain different in 
$\Lambda^{\mathrm{st}}_{\beta*\gamma}$ since the total number of such blocks 
equals  $l$. Thus
$\Lambda^{\mathrm{st}}_{\beta}=\Lambda^{\mathrm{st}}_{\beta*\gamma}$.
We also have  $\ran^{\mathrm{st}}(\beta*\gamma)=\ran^{\mathrm{st}}(\gamma)$ 
as both sets have the same cardinality. This completes the proof.
\end{proof}

Let $\mathcal{X}$ be the set of all unordered partitions of $N$ into $l$ disjoint 
blocks satisfying the condition that all $a_i$'s belong to different blocks. We 
interpret the elements of $\mathcal{X}$ as equivalence relations on $N$. Let further
$\mathcal{Y}$ be the set of all sets of the form $\{b_1,\dots,b_l\}$ where
$b_i\in A_i$ for all $i=1,\dots,l$. Now we are ready to describe the minimal
isolated subsemigroups in $\mathcal{T}_n^{(l)}$.

\begin{proposition}\label{pr12}
\begin{enumerate}[(a)]
\item\label{pr12.1} For every $\Lambda\in \mathcal{X}$ and $I\in \mathcal{Y}$ 
there exists a unique idempotent, $\varepsilon(\Lambda,I)\in (\mathcal{T}_n,*)$, of 
stable rank $l$ such that $\Lambda_{\varepsilon(\Lambda,I)}=\Lambda$ and
$\ran(\varepsilon(\Lambda,I))=I$.
\item\label{pr12.2} For every $\Lambda\in \mathcal{X}$ and $I\in \mathcal{Y}$
the set $\sqrt{\varepsilon(\Lambda,I)}$ is an isolated subsemigroup of 
$(\mathcal{T}_n,*)$ and $\sqrt{\varepsilon(\Lambda,I)}=\mathtt{I}(\beta)$
for every $\beta\in \sqrt{\varepsilon(\Lambda,I)}$.
\end{enumerate}
\end{proposition}

\begin{proof}
The statement \eqref{pr12.1} follows from the description of idempotents in 
$(\mathcal{T}_n,*)$, obtained in Theorem~\ref{t6}. 

To prove \eqref{pr12.2} we first note that, by definition, the set 
$\sqrt{\varepsilon(\Lambda,I)}$ contains exactly elements with stable kernel 
$\Lambda$ and stable image $I$. Let 
$\beta,\gamma\in \sqrt{\varepsilon(\Lambda,I)}$.
Then Lemma~\ref{l11}\eqref{l11.2} implies that $\beta*\gamma\in 
\sqrt{\varepsilon(\Lambda,I)}$ as well, in particular, $\sqrt{\varepsilon(\Lambda,I)}$
is a subsemigroup of $(\mathcal{T}_n,*)$. Hence \eqref{pr12.2} follows from
Proposition~\ref{pr4}. This completes the proof.
\end{proof}

\begin{proposition}\label{pr14}
\begin{enumerate}[(a)]
\item\label{pr14.1} $\mathcal{T}_n^{(l)}$ is a subsemigroup of $(\mathcal{T}_n,*)$. 
\item\label{pr14.2} 
The decomposition 
\begin{equation}\label{eq.pr14.2}
\mathcal{T}_n^{(l)}=
\bigcup_{(\Lambda,I)\in\mathcal{X}\times\mathcal{Y}}
\sqrt{\varepsilon(\Lambda,I)}
\end{equation}
into a disjoint union of subsets defines a congruence on 
$\mathcal{T}_n^{(l)}$, the quotient modulo which is isomorphic 
to the rectangular band $\mathcal{X}\times\mathcal{Y}$. 
\item\label{pr14.3} Let $\psi:\mathcal{T}_n^{(l)}\to \mathcal{X}\times\mathcal{Y}$ be 
the natural projection, given by \eqref{pr14.2}. Then
$\psi^{-1}$ induces a bijection between the set of all subsemigroups of
$\mathcal{X}\times\mathcal{Y}$ and $\mathfrak{I}(\mathcal{T}_n^{(l)})$.
\item\label{pr14.4} $|\mathfrak{I}(\mathcal{T}_n^{(l)})|=2^{|\mathcal{X}|+|\mathcal{Y}|}-
2^{|\mathcal{X}|}-2^{|\mathcal{Y}|}+1$.
\end{enumerate}
\end{proposition}

We remark that $|\mathcal{X}|=l^{n-l}$ and $|\mathcal{Y}|=\prod_{i=1}^l|A_i|$.
Note also that $T\subset \mathcal{X}\times\mathcal{Y}$ is a
subsemigroup if and only if the following condition is satisfied:
for all $(\Lambda,I),(\Lambda',I')\in \mathcal{X}\times\mathcal{Y}$
we have that $(\Lambda,I),(\Lambda',I')\in T$ implies both
$(\Lambda,I')\in T$ and $(\Lambda',I)\in T$.

\begin{proof}
The statements \eqref{pr14.1} and \eqref{pr14.2} follow from 
Lemma~\ref{l11}\eqref{l11.2}.
The rectangular band $\mathcal{X}\times\mathcal{Y}$ consists of idempotents,
and hence every subsemigroup of it is isolated. Now \eqref{br8} implies that
$\psi^{-1}$ (i.e. taking the pre-image with respect to  $\psi$) 
induces a map from  $\mathfrak{I}(\mathcal{X}\times\mathcal{Y})$
to $\mathfrak{I}(\mathcal{T}_n^{(l)})$, which is injective since $\psi$ is
surjective. Let now $T\in \mathfrak{I}(\mathcal{T}_n^{(l)})$. Since the congruence
classes are minimal isolated subsemigroups of $\mathcal{T}_n^{(l)}$ by 
Proposition~\ref{pr12}, it follows that $T$ must be the union of those congruence
classes, which have a non-trivial intersection with $T$. Hence $\psi(T)$ is 
a subsemigroup of $\mathcal{X}\times\mathcal{Y}$ and $\psi^{-1}(\psi(T))=T$, 
implying that the map, induced by $\psi^{-1}$, is even bijective. 
This proves \eqref{pr14.3}.

Because of  \eqref{pr14.3} we just compute the cardinality of the set
$\mathfrak{T}(\mathcal{X}\times\mathcal{Y})$ of all subsemigroups in
$\mathcal{X}\times\mathcal{Y}$. It is easy to see that every subsemigroup in
$\mathcal{X}\times\mathcal{Y}$ is uniquely determined by the sets of the
left and the right coordinates of its elements. So, there are exactly
$(2^{|\mathcal{X}|}-1)(2^{|\mathcal{Y}|}-1)$ subsemigroups
in $\mathcal{X}\times\mathcal{Y}$. The statement \eqref{pr14.4} follows and
the proof of our proposition is complete.
\end{proof}

For $\varnothing\neq X\subset \mathcal{X}$ and $\varnothing\neq 
Y\subset\mathcal{Y}$ let $F(X,Y)$ denote the union of all 
$\sqrt{\varepsilon}$ for which $\Lambda_{\varepsilon}\in X$ and 
$\im(\varepsilon)\in Y$. From Proposition~\ref{pr14} we obtain:

\begin{corollary}\label{ccc123}
Every isolated subsemigroup of $(\mathcal{T}_n,*)$, which is contained
in $\mathcal{T}_n^{(l)}$, is equal to $F(X,Y)$ for appropriate $X$ and $Y$.
\end{corollary}

\subsection{The case of stable rank $l-1$}\label{s4.2}

Let $\varepsilon\in E(\mathcal{T}_n,*)$ be some idempotent of rank $l-1$,
given by \eqref{idemp}. From Theorem~\ref{t6} it follows that 
there exists a unique element, $i\in\{1,\dots,l-1\}$, and a unique pair,
$\{m,k\}\subset\{1,\dots,l\}$, $m\neq k$, such that 
$\{a_m,a_k\}\in E_i$ and $e_i\in A_m$.
We call the pair $\{a_m,a_k\}$ {\em distinguished}, the element
$a_k$ {\em trifle} and the element $a_m$ {\em burdened} (for the idempotent
$\varepsilon$) and denote them by $\mathfrak{d}(\varepsilon)$, 
$\mathfrak{t}(\varepsilon)$ and $\mathfrak{b}(\varepsilon)$ respectively. 
By $\mathfrak{nt}(\varepsilon)$ and $\mathfrak{nb}(\varepsilon)$ we denote 
the elements of $\{1,\dots,l\}$ such that 
$\mathfrak{t}(\varepsilon)=a_{\mathfrak{nt}(\varepsilon)}$ and
$\mathfrak{b}(\varepsilon)=a_{\mathfrak{nb}(\varepsilon)}$ respectively.
It will be convenient now to adjust \eqref{idemp} in the following way:
if necessary, we can rearrange the indexes of $e_i$'s to obtain 
\begin{equation}\label{idempnew}
\varepsilon=\left(
\begin{array}{cccccc}
E_1 & \dots &  E_{\mathfrak{nt}(\varepsilon)-1} & 
E_{\mathfrak{nt}(\varepsilon)+1} & \dots & E_l \\
e_1 & \dots &  e_{\mathfrak{nt}(\varepsilon)-1} & 
e_{\mathfrak{nt}(\varepsilon)+1} & \dots & e_l \\
\end{array}
\right),
\end{equation}
where $a_i\in E_i$ for all $i=1,\dots,\mathfrak{nt}(\varepsilon)-1,
\mathfrak{nt}(\varepsilon)+1,\dots,l$. Then, in particular,
$\mathfrak{t}(\varepsilon)\in E_{\mathfrak{nb}(\varepsilon)}$.

For $\beta\in\sqrt{\varepsilon}$ we set 
$\mathfrak{d}(\beta)=\mathfrak{d}(\varepsilon)$,
$\mathfrak{t}(\beta)=\mathfrak{t}(\varepsilon)$,
$\mathfrak{b}(\beta)=\mathfrak{b}(\varepsilon)$,
$\mathfrak{nt}(\beta)=\mathfrak{nt}(\varepsilon)$, and
$\mathfrak{nb}(\beta)=\mathfrak{nb}(\varepsilon)$. 
For $k=1,\dots,l$ let $\mathcal{T}_n^{(l-1,k)}$ denote the set of all 
elements $\beta\in \mathcal{T}_n^{(l-1)}$ such that $\mathfrak{t}(\beta)=a_k$.
For $k,m=1,\dots,l$, $k\neq m$, let $\mathcal{T}_n^{(l-1,k,m)}$ denote the set 
of all elements $\beta\in \mathcal{T}_n^{(l-1,k)}$ such that 
$\mathfrak{b}(\beta)=a_m$.

\begin{lemma}\label{l20}
Let $\varepsilon$ be an idempotent in $(\mathcal{T}_n,*)$ of
rank $l-1$ given by \eqref{idempnew} and $\beta\in\sqrt{\varepsilon}$. Then
\begin{enumerate}[(a)]
\item\label{l20.1}
$\im^{\mathrm{st}}(\beta)=\beta(A\setminus\{\mathfrak{t}(\beta)\})$.
\item\label{l20.2} If $\rank(\beta\alpha)=l-1$ then
$\Lambda_{\varepsilon}=\Lambda_{\beta}^{\mathrm{st}}=\Lambda_{\beta\alpha}$.
\item\label{l20.3} If $\rank(\beta\alpha)=l$ then
every block of $\Lambda_{\beta\alpha}$ contains at most two
elements from $A$, moreover, there is the unique block, namely
$(\beta\alpha)^{-1}(\mathfrak{t}(\beta))$, which does not contain
any element from $A$, and there is the unique block which contains
exactly two elements from  $A$, namely the elements from
$\mathfrak{d}(\beta)$.
\item\label{l20.4} $\beta\alpha$ induces a permutation on
$A\setminus\{\mathfrak{t}(\beta)\}$
\end{enumerate}
\end{lemma}

\begin{proof}
We obviously have $\Lambda_{\beta\alpha}\subset \Lambda_{\beta}^{\mathrm{st}}$
and in the case $\rank(\beta\alpha)=l-1$ the equality follows since the number
of blocks in $\Lambda_{\beta\alpha}$ and $\Lambda_{\beta}^{\mathrm{st}}$ 
coincide. This proves \eqref{l20.2}. Moreover, in this case
$\im^{\mathrm{st}}(\beta)=\beta(A\setminus\{\mathfrak{t}(\beta)\})$ follows 
immediately from the definition of $\mathfrak{t}(\beta)$, which proves the
corresponding part of \eqref{l20.1}.

Since $\str(\beta)\neq l$, in the case  $\rank(\beta\alpha)=l$ it is not 
possible that every block of $\Lambda_{\beta\alpha}$ contains an element from 
$A$. On the other hand, every block of $\Lambda_{\beta}^{\mathrm{st}}$ 
contains at least one element from $A$ and hence, since $\str(\beta)= l-1$, 
it follows that all blocks of $\Lambda_{\beta}^{\mathrm{st}}$ except one 
contain exactly one element from $A$, and the exceptional
block contains two elements from $A$. This implies the uniqueness of
the block with two elements from $A$ and of the block without any elements from
$A$ for $\Lambda_{\beta\alpha}$. From the definition of $\mathfrak{t}(\beta)$
it follows that $(\beta\alpha)^{-1}(\mathfrak{t}(\beta))$ 
does not contain any element of $A$, and from the
definition of $\mathfrak{d}(\beta)$ it follows that $\mathfrak{d}(\beta)$ is
contained in some block of $\Lambda_{\beta\alpha}$. This proves \eqref{l20.3}. 
Since the pre-image of $\mathfrak{t}(\beta)$ does not intersect $A$, 
$\mathfrak{t}(\beta)$ can not belong to $\im^{\mathrm{st}}(\beta\alpha)$, 
which completes the proof of \eqref{l20.1} as well. Finally, \eqref{l20.4} follows 
from \eqref{l20.1} and the lemma is proved.
\end{proof}

\begin{corollary}\label{c201}
Let $\beta,\gamma\in\mathcal{T}_n^{(l-1)}$. Then 
$\rank(\beta\alpha\gamma\alpha)\geq \rank(\beta\alpha)-1$ and
$\Lambda_{\beta\alpha\gamma\alpha}$ is obtained from
$\Lambda_{\beta\alpha}$ by uniting $(\beta\alpha)^{-1}(\mathfrak{b}(\gamma))$
and $(\beta\alpha)^{-1}(\mathfrak{t}(\gamma))$.
\end{corollary}

\begin{proof}
$\beta\alpha$ maps the blocks of $\Lambda_{\beta\alpha}$ to
$\im(\beta\alpha)\subset A$. Computing $\Lambda_{\beta\alpha\gamma\alpha}$
we should unite those blocks of $\Lambda_{\beta\alpha}$, whose images belong
to the same block of $\Lambda_{\gamma\alpha}$. However,
Lemma~\ref{l20}\eqref{l20.3} says that there is a unique block in
$\Lambda_{\gamma\alpha}$, which contains more than one element from $A$,
and this block contains $\mathfrak{d}(\gamma)$. The statement follows.
\end{proof}

\begin{corollary}\label{c202}
Let $\beta,\gamma\in\mathcal{T}_n^{(l-1)}$ be such that 
$\beta*\gamma\in\mathcal{T}_n^{(l-1)}$. Then
\begin{enumerate}[(a)]
\item\label{c202.1} $\mathfrak{d}(\beta*\gamma)=\mathfrak{d}(\beta)$.
\item\label{c202.2} $\mathfrak{t}(\beta*\gamma)=\mathfrak{t}(\gamma)$.
\end{enumerate}
\end{corollary}

\begin{proof}
\eqref{c202.1} follows from Lemma~\ref{l20}\eqref{l20.3}.
To prove \eqref{c202.2} consider $\beta\alpha\gamma\alpha$.
From Lemma~\ref{l20}\eqref{l20.3} we have that
$(\gamma\alpha)^{-1}(\mathfrak{t}(\gamma))\cap A=\varnothing$, which obviously
implies that $(\beta\alpha\gamma\alpha)^{-1}(\mathfrak{t}(\gamma))=\varnothing$. 
The statement now follows from the definition of the trifle element.
\end{proof}

Now we can describe all minimal isolated subsemigroups in $\mathcal{T}_n^{(l-1)}$.

\begin{corollary}\label{c21}
Let $\varepsilon\in \mathcal{T}_n^{(l-1)}$ be  an idempotent given by 
\eqref{idempnew}. 
Then $\sqrt{\varepsilon}$ is a minimal isolated subsemigroup
of $(\mathcal{T}_n,*)$.
\end{corollary}

\begin{proof}
Because of Proposition~\ref{pr4} we need only to prove that
$\sqrt{\varepsilon}$ is a semigroup. 
Let $\beta,\gamma\in \sqrt{\varepsilon}$. 
We obviously have $\Lambda_{\beta\alpha}\subset
\Lambda_{\beta\alpha\gamma\alpha}$. By Corollary~\ref{c201},
$\Lambda_{\beta\alpha\gamma\alpha}$ is obtained from 
$\Lambda_{\beta\alpha}$ uniting the pre-images of
$\mathfrak{t}(\gamma)$($=\mathfrak{t}(\beta)$) and
$\mathfrak{b}(\gamma)$($=\mathfrak{b}(\beta)$).
By Lemma~\ref{l20}\eqref{l20.4}, the element 
$\beta\alpha\gamma\alpha$ induces a permutation on 
$A\setminus \{\mathfrak{t}(\beta)\}$,
in particular, $\rank(\beta\alpha\gamma\alpha)\geq l-1$. On the
other hand, as $\mathfrak{d}(\beta)=\mathfrak{d}(\gamma)$, it follows
that $\rank(\beta\alpha\gamma\alpha)\leq l-1$. Hence
$\rank(\beta\alpha\gamma\alpha)= l-1$. Note
that $(\beta\alpha\gamma\alpha)^{*m}=(\beta\alpha\gamma\alpha)^m$
since $\alpha\alpha=\alpha$.
Applying the previous arguments inductively we get 
$\Lambda_{(\beta\alpha\gamma\alpha)^m}=\Lambda_{\beta\alpha\gamma\alpha}$
and $\im((\beta\alpha\gamma\alpha)^m)=\im(\beta\alpha\gamma\alpha)$
for all $m$, in particular, it follows that the stable rank of
$\beta\alpha\gamma\alpha$ equals $l-1$. As both $\beta\alpha\gamma\alpha$ 
and $\gamma$ have stable rank $l-1$ we also deduce that 
$\im^{\mathrm{st}}(\beta*\gamma)=\im^{\mathrm{st}}(\gamma)$.
Using $\mathfrak{d}(\beta)=\mathfrak{d}(\gamma)$, we
further get
$\Lambda_{\beta*\gamma}^{\mathrm{st}}=\Lambda_{\beta}^{\mathrm{st}}=
\Lambda_{\beta\alpha\gamma\alpha}$.
In particular, if $m$ is such that $(\beta\alpha\gamma\alpha)^{*m}=\nu$
is an idempotent, then $\Lambda_{\nu}=\Lambda_{\varepsilon}$ and
$\im^{\mathrm{st}}(\nu)=\im^{\mathrm{st}}(\varepsilon)$, which implies
that $\nu=\varepsilon$ since every idempotent in 
$(\mathcal{T}_n,*)$ is uniquely defined by the corresponding $\Lambda$
and the image, see Theorem~\ref{t6}. Therefore $\beta*\gamma\in \sqrt{\varepsilon}$
and hence  $\sqrt{\varepsilon}$ is indeed a semigroup. The statement follows.
\end{proof}

\begin{proposition}\label{p18}
\begin{enumerate}[(a)]
\item\label{p18.1} For every $k,m\in\{1,\dots,l\}$, $k\neq m$, the set 
$\mathcal{T}_n^{(l-1,k,m)}$ is an isolated subsemigroup of $(\mathcal{T}_n,*)$.
\item\label{p18.2} For every $k\in\{1,\dots,l\}$ the set $\mathcal{T}_n^{(l-1,k)}$ is
an isolated subsemigroup of $(\mathcal{T}_n,*)$.
\item\label{p18.3} For every $k,m\in\{1,\dots,l\}$, $k\neq m$, the set 
$\mathcal{T}_n^{(l-1,k,m)}\cup
\mathcal{T}_n^{(l-1,m,k)}$ is an isolated subsemigroup of $(\mathcal{T}_n,*)$.
\item\label{p18.4} Let $\varepsilon_1\neq \varepsilon_2$ be two
idempotents in $\mathcal{T}_n^{(l-1)}$. Assume that
$\mathfrak{d}(\varepsilon_1)\neq\mathfrak{d}(\varepsilon_2)$ 
and $\mathfrak{t}(\varepsilon_1)\neq\mathfrak{t}(\varepsilon_2)$.
Then $\sqrt{\varepsilon_1}*\sqrt{\varepsilon_2}$ is not contained
in $\mathcal{T}_n^{(l-1)}$.
\end{enumerate}
\end{proposition}

\begin{proof}
The statement \eqref{p18.1} follows from \eqref{p18.2} and \eqref{p18.3} using \eqref{br1}.

To prove \eqref{p18.2} let us first consider the case when
$\alpha=\mathrm{id}$. In this case the only idempotents $\mathcal{T}_n^{(l-1,k)}$
contains are $\varepsilon_{m,k}$, $m\in\{1,\dots,k-1,k+1,\dots,n\}$,
defined as follows:
\begin{equation}\label{p18.eq.1}
\varepsilon_{m,k}=\left(
\begin{array}{cccccccc}
1 & 2 & \dots & k-1 & k & k+1 & \dots & n \\
1 & 2 & \dots & k-1 & m & k+1 & \dots & n
\end{array}
\right).
\end{equation}
If $\beta\in \mathcal{T}_n$ is such that
$\beta^j=\varepsilon_{m,k}$ for some $m$ and $j>0$ then 
$\rank(\beta)=n-1$ and $\beta$ induces a permutation on
$N\setminus\{k\}$. On the other hand, if $\beta$ is not invertible and
induces a permutation on $N\setminus\{k\}$, it must have rank $n-1$.
It follows that the set $\cup_{m\neq k}\sqrt{\varepsilon_{m,k}}$
consists of all noninvertible elements in $\mathcal{T}_n$, which 
induce a permutation on $N\setminus\{k\}$. In particular, this set is
closed with respect to composition and hence is a subsemigroup.
Since it is a union of isolated subsemigroups by Corollary~\ref{c21}, 
it is an isolated subsemigroup itself
by Proposition~\ref{pr3}. Hence \eqref{p18.2} is true in the case of invertible
$\alpha$. In the general case consider the homomorphism $\overline{\varphi}:
(\mathcal{T}_n,*)\to \mathcal{T}(A)$. Let $\mathcal{T}(A)^{(l-1,k)}$ denote
the set of all elements in $\mathcal{T}(A)$ of stable rank $l-1$ and with
the trifle element $a_k$. Let $\beta\in \mathcal{T}_n^{(l-1,k)}$. Then
from Lemma~\ref{l20} it follows that $\alpha\beta\alpha$ has rank
$l-1$ and induces a permutation on $A\setminus\{a_k\}$. In particular,
$\overline{\varphi}(\beta)\in \mathcal{T}(A)^{(l-1,k)}$. On the other
hand, if $\overline{\varphi}(\beta)\in \mathcal{T}(A)^{(l-1,k)}$
then $\str(\beta)=l-1$ and $\mathrm{t}(\beta)=a_k$, that is
$\beta\in \mathcal{T}_n^{(l-1,k)}$. This implies that
$\mathcal{T}_n^{(l-1,k)}=\overline{\varphi}^{-1}(\mathcal{T}(A)^{(l-1,k)})$.
Since $\mathcal{T}(A)^{(l-1,k)}$ is an isolated subsemigroup of $\mathcal{T}(A)$,
we get, using \eqref{br8}, that $\mathcal{T}_n^{(l-1,k)}$ is an isolated subsemigroup 
of $(\mathcal{T}_n,*)$. This proves \eqref{p18.2}. The statement \eqref{p18.3} 
is proved by analogous arguments (that is one first checks it in the case
$\alpha=\mathrm{id}$ and then uses $\overline{\varphi}^{-1}$).

Finally, to prove  \eqref{p18.4} we show that
$\str(\varepsilon_1*\varepsilon_2)<l-1$ if 
$\mathfrak{d}(\varepsilon_1)\neq\mathfrak{d}(\varepsilon_2)$ 
and $\mathfrak{t}(\varepsilon_1)\neq\mathfrak{t}(\varepsilon_2)$.
Obviously $\str(\varepsilon_1*\varepsilon_2)\leq
\rank(\overline{\varphi}(\varepsilon_1)\overline{\varphi}(\varepsilon_2))$ and
$\overline{\varphi}(\varepsilon_i)=
\varepsilon_{\mathfrak{b}(\varepsilon_i),\mathfrak{t}(\varepsilon_i)}$ 
for $i=1,2$ (here we use the notation similar to \eqref{p18.eq.1}). 
Under our assumption a direct calculation shows that
$\rank((\varepsilon_{\mathfrak{b}(\varepsilon_1),\mathfrak{t}(\varepsilon_1)}
\varepsilon_{\mathfrak{b}(\varepsilon_2),\mathfrak{t}(\varepsilon_2)})^2)=l-2$.
The statement \eqref{p18.4} follows and hence the proof is complete.
\end{proof}

Now we would like to present various constructions of isolated subsemigroups 
in $\mathcal{T}_n^{(l-1)}$. We start with the one which corresponds to
the case described in Proposition~\ref{p18}\eqref{p18.1}.

Fix $k,m\in\{1,\dots,l\}$, $k\neq m$. Let $\mathcal{U}^{(k,m)}$ denote the set 
of all unordered partitions of $N$ into $l-1$ disjoint blocks satisfying the
condition that every block contains at least one element from $A$ and that
$a_k$ and $a_m$ are contained in the same block. Note that $\mathcal{U}^{(k,m)}=
\mathcal{U}^{(m,k)}$. Let further
$\mathcal{V}^{(k)}$ denote the set of all sets of the form
$\{b_1,\dots,b_{k-1},b_{k+1},\dots,b_l\}$, where $b_i\in A_i$ for all
$i=1,\dots,k-1,k+1,\dots,l$.

For $k,m\in\{1,\dots,l\}$, $k\neq m$, $\varnothing\neq X\subset \mathcal{U}^{(k,m)}$,
and $\varnothing\neq Y\subset\mathcal{V}^{(k)}$, let $H(k,m,X,Y)$ denote the union of all
$\sqrt{\varepsilon}$, where $\varepsilon\in \mathcal{T}_n^{(l-1,k,m)}$ is
an idempotent such that $\Lambda_{\varepsilon}\in X$ and
$\im(\varepsilon)\in Y$. For $H(k,m,X,Y)$ we can prove an analogue of Proposition~\ref{pr14}
and Corollary~\ref{ccc123}.

\begin{proposition}\label{pr301}
Let $k,m\in\{1,\dots,l\}$, $k\neq m$.
\begin{enumerate}[(a)]
\item\label{pr301.1} For every $\Lambda\in \mathcal{U}^{(k,m)}$ and $I\in \mathcal{V}^{(k)}$
there exists a unique idempotent, $\varepsilon(\Lambda,I)$, such that
$\Lambda_{\varepsilon(\Lambda,I)}=\Lambda$ and $\im(\varepsilon(\Lambda,I))=I$.
\item\label{pr301.2} 
The decomposition
\begin{displaymath}
\mathcal{T}_n^{(l-1,k,m)}=\bigcup_{(\Lambda,I)\in\mathcal{U}^{(k,m)}\times 
\mathcal{V}^{(k)}}\sqrt{\varepsilon(\Lambda,I)}
\end{displaymath}
into a disjoint union of subsets defines a congruence on
$\mathcal{T}_n^{(l-1,k,m)}$, the quotient modulo which is isomorphic
to the rectangular band $\mathcal{U}^{(k,m)}\times 
\mathcal{V}^{(k)}$.
\item\label{pr301.3}
For every $\varnothing\neq X\subset \mathcal{U}^{(k,m)}$ and 
$\varnothing\neq Y\subset \mathcal{V}^{(k)}$ the set 
$H(k,m,X,Y)$ is an isolated subsemigroup of $(\mathcal{T}_n,*)$.
\item\label{pr301.4}
Every isolated subsemigroup of $\mathcal{T}_n^{(l-1,k,m)}$ is equal to
$H(k,m,X,Y)$ for appropriate $X$ and $Y$.
\item\label{pr301.5}
$\mathcal{T}_n^{(l-1,k,m)}$ contains 
$2^{|\mathcal{U}^{(k,m)}|+|\mathcal{V}^{(k)}|}-
2^{|\mathcal{U}^{(k,m)}|}-2^{|\mathcal{V}^{(k)}|}+1$ isolated subsemigroups.
\end{enumerate}
\end{proposition}

\begin{proof}
The statement \eqref{pr301.1} follows from Theorem~\ref{t6}. Let 
$(\Lambda,I),(\Lambda',I')\in \mathcal{U}^{(k,m)}\times \mathcal{V}^{(k)}$,
$\beta\in \sqrt{\varepsilon(\Lambda,I)}$ and $\gamma\in \sqrt{\varepsilon(\Lambda',I')}$.
From Proposition~\ref{p18}\eqref{p18.1} we have $\beta*\gamma\in \mathcal{T}_n^{(l-1,k,m)}$,
in particular, $\str(\beta*\gamma)=l-1$. Hence, using 
$\mathfrak{b}(\beta)=\mathfrak{b}(\gamma)$, $\mathfrak{t}(\beta)=\mathfrak{t}(\gamma)$, 
and  Corollary~\ref{c201}, we get $\Lambda^{\mathrm{st}}_{\beta*\gamma}=
\Lambda^{\mathrm{st}}_{\beta}$. Moreover, Corollary~\ref{c202}\eqref{c202.2}
and Lemma~\ref{l20}\eqref{l20.1} imply $\im^{\mathrm{st}}(\beta*\gamma)=
\im^{\mathrm{st}}(\gamma)$. Thus $\beta*\gamma\in \sqrt{\varepsilon(\Lambda,I')}$,
which implies \eqref{pr301.2}. The rest is proved by the same arguments as in
Proposition~\ref{pr14} and Corollary~\ref{ccc123}.
\end{proof}

Our next construction corresponds to
the case described in Proposition~\ref{p18}\eqref{p18.3}.
For $k,m\in\{1,\dots,l\}$, $k\neq m$, and   $\varnothing\neq X\subset \mathcal{U}^{(k,m)}=\mathcal{U}^{(m,k)}$, let $K(\{k,m\},X)$ denote 
the union of all $\sqrt{\varepsilon}$, where 
$\varepsilon\in  \mathcal{T}_n^{(l-1,k,m)}\cup 
\mathcal{T}_n^{(l-1,m,k)}$ is an idempotent such that 
$\Lambda_{\varepsilon}\in X$.

\begin{proposition}\label{pr302}
Let $k,m\in\{1,\dots,l\}$, $k\neq m$.
\begin{enumerate}[(a)]
\item\label{pr302.1}
For every $X$ as above the set $K(\{k,m\},X)$ is an 
isolated subsemigroup of $\mathcal{T}_n^{(l-1,k,m)}\cup \mathcal{T}_n^{(l-1,m,k)}$.
\item\label{pr302.2}
Every isolated subsemigroup of $\mathcal{T}_n^{(l-1,k,m)}\cup \mathcal{T}_n^{(l-1,m,k)}$,
which intersects both $\mathcal{T}_n^{(l-1,k,m)}$ and $\mathcal{T}_n^{(l-1,m,k)}$ in a
non-trivial way, is equal to the semigroup 
$K(\{k,m\},X)$ for an appropriate $X$.
\item\label{pr302.3}
The number of isolated subsemigroups of 
$\mathcal{T}_n^{(l-1,k,m)}\cup 
\mathcal{T}_n^{(l-1,m,k)}$, which intersect both subsemigroups
$\mathcal{T}_n^{(l-1,k,m)}$ and $\mathcal{T}_n^{(l-1,m,k)}$ in a non-trivial way, equals $2^{\displaystyle (l-1)^{n-l}}-1$.
\end{enumerate}
\end{proposition}

\begin{proof}
That both $K(\{k,m\},X)\cap \mathcal{T}_n^{(l-1,k,m)}$ and
$K(\{k,m\},X)\cap \mathcal{T}_n^{(l-1,m,k)}$ are subsemigroups
follows from Proposition~\ref{pr301}. Let $\Lambda,\Lambda'\in X$, 
$I\in \mathcal{V}^{(k)}$, $I'\in \mathcal{V}^{(m)}$, 
$\beta\in \sqrt{\varepsilon(\Lambda,I)}$ and 
$\gamma\in \sqrt{\varepsilon(\Lambda',I')}$. 
From Proposition~\ref{p18}\eqref{p18.3} 
we have $\str(\beta*\gamma)=l-1$. Hence, using 
$\mathfrak{d}(\beta)=\mathfrak{d}(\gamma)$ and  Corollary~\ref{c201}, 
we get $\Lambda^{\mathrm{st}}_{\beta*\gamma}=
\Lambda^{\mathrm{st}}_{\beta}=\Lambda\in X$.  Further,
$\mathfrak{d}(\beta*\gamma)=\mathfrak{d}(\beta)$
by Corollary~\ref{c202}\eqref{c202.1}. Thus
$\beta*\gamma\in K(\{k,m\},X)$ and the claim \eqref{pr302.1} follows.

Let now $T$ be an isolated subsemigroup of $\mathcal{T}_n^{(l-1,k,m)}\cup 
\mathcal{T}_n^{(l-1,m,k)}$, which intersects both $\mathcal{T}_n^{(l-1,k,m)}$ and $\mathcal{T}_n^{(l-1,m,k)}$ in a non-trivial way. Set
\begin{displaymath}
X=\{\Lambda^{\mathrm{st}}_{\beta}:\beta\in T\}.
\end{displaymath}
Then $T\subset K(\{k,m\},X)$. 

Let now $\Lambda\in X$ and $\beta\in T$ be an idempotent such that 
$\Lambda_{\beta}=\Lambda$. Assume that $\beta\in \mathcal{T}_n^{(l-1,k,m)}$.
Then for any $\gamma\in \mathcal{T}_n^{(l-1,m,k)}\cap T$ from 
Corollary~\ref{c202} it follows that $\beta*\gamma\in 
\mathcal{T}_n^{(l-1,m,k)}$, and as both $\beta$ and $\beta*\gamma$ have
rank $l-1$, it also follows that $\Lambda_{\beta}=\Lambda_{\beta*\gamma}$.
This means that
\begin{equation}\label{eqnewe1}
X=\{\Lambda^{\mathrm{st}}_{\beta}:\beta\in T\cap
\mathcal{T}_n^{(l-1,m,k)}\}=
\{\Lambda^{\mathrm{st}}_{\beta}:\beta\in T\cap
\mathcal{T}_n^{(l-1,k,m)}\}.
\end{equation}

Now let $\varepsilon$ be 
an idempotent of $\mathcal{T}_n^{(l-1,k,m)}\cap T$ given by \eqref{idempnew}. 
Let $i\in\{1,\dots,\mathfrak{nt}(\varepsilon)-1,\mathfrak{nt}(\varepsilon)+1,
\dots,l\}$ and $e'_i\in E_i$ be some element different from $e_i$. Consider
the idempotent $\varepsilon'$ defined as follows:
\begin{displaymath}
\varepsilon'(x)=
\begin{cases}
\varepsilon(x),& x\not\in E_i;\\
e'_i, & \text{otherwise.}
\end{cases}
\end{displaymath}
Obviously $\Lambda_{\varepsilon'}=\Lambda_{\varepsilon}$.

\begin{lemma}\label{lnew-l1}
$\varepsilon'\in T$. 
\end{lemma}

\begin{proof}
Because of \eqref{eqnewe1} we can fix an idempotent, 
$\beta\in \mathcal{T}_n^{(l-1,m,k)}\cap T$, 
such that $\Lambda_{\beta}=\Lambda_{\varepsilon}$. 
Consider now the element $\gamma\in\mathcal{T}_n$ defined as follows:
\begin{displaymath}
\gamma(x)=
\begin{cases}
e'_i, & x=\mathfrak{t}(\varepsilon);\\
e_i, & x\in E_{\mathfrak{nb}(\varepsilon)},
x\neq\mathfrak{t}(\varepsilon);\\
e_{\mathfrak{nb}(\varepsilon)}, & x\in E_i,
x\neq\mathfrak{t}(\varepsilon);\\
\varepsilon(x),& \text{ otherwise.}
\end{cases}
\end{displaymath}
From the definitions we have $\gamma*\gamma=\varepsilon$ and hence
$\gamma\in T$. In particular, $\beta*\gamma\in T$. As both
$\beta*\gamma$ and $\beta$ have rank $l-1$, we have
$\Lambda_{\beta*\gamma}=\Lambda_{\beta}=
\Lambda_{\varepsilon}=\Lambda_{\varepsilon'}$. We further have:
\begin{displaymath}
\begin{array}{lcl}
\im^{\mathrm{st}}(\beta*\gamma)&=&\\
\text{(by Lemma~\ref{l20}\eqref{l20.1})}
&=& \beta*\gamma(A\setminus
\{\mathfrak{t}(\beta*\gamma)\})\\
\text{(by Corollary~\ref{c202}\eqref{c202.2})}&=&\beta*\gamma(A\setminus
\{\mathfrak{t}(\gamma)\})\\
\text{(by our choice of $\beta$ and $\gamma$)}&=&\beta*\gamma(A\setminus
\{\mathfrak{b}(\beta)\})\\
\text{(by Lemma~\ref{l20})}&=&\beta\alpha\gamma(A)\\
\text{(by the definition of $\mathfrak{t}(\beta)$)}
&=&\gamma(A\setminus\{\mathfrak{t}(\beta)\})\\
\text{(by our choice of $\beta$ and $\gamma$)}
&=&\gamma(A\setminus\{\mathfrak{b}(\gamma)\})\\
\text{(by Lemma~\ref{l20})}
&=&\big(\gamma(A\setminus\{\mathfrak{t}(\gamma)\})\setminus
\gamma(\{\mathfrak{b}(\gamma)\})\big)\cup \gamma(\{\mathfrak{t}(\gamma)\})\\
\text{(by Lemma~\ref{l20}\eqref{l20.1})}
&=&\big(\im^{\mathrm{st}}(\gamma)\setminus
\gamma(\{\mathfrak{b}(\gamma)\})\big)\cup \gamma(\{\mathfrak{t}(\gamma)\})\\
\text{(by the definition of $\gamma$)}
&=&\big(\im^{\mathrm{st}}(\gamma)\setminus\{e_i\}\big)\cup \{e'_i\}\\
\text{(by the definition of $\varepsilon'$)}
&=&\im^{\mathrm{st}}(\varepsilon').
\end{array}
\end{displaymath}

Thus the element $\beta*\gamma\in T$ satisfies 
$\Lambda_{\beta*\gamma}=\Lambda_{\varepsilon'}$ and
$\im^{\mathrm{st}}(\beta*\gamma)=\im^{\mathrm{st}}(\varepsilon')$
implying $\beta*\gamma\in\sqrt{\varepsilon'}$, which means 
$\varepsilon'\in T$.
\end{proof}

Lemma~\ref{lnew-l1} basically says that if $T$ contains an idempotent with
a certain stable kernel and a certain stable image, the $T$ also 
contains an idempotent with the same stable kernel, but in which the 
stable image is different by one element,
arbitrarily chosen in arbitrary block $E_i$. From this and the definition 
of the semigroup  $K(\{k,m\},X)$ it follows easily  that 
$\mathcal{T}_n^{(l-1,k,m)}\cap K(\{k,m\},X)\subset T$. Analogously
we obtain 
$\mathcal{T}_n^{(l-1,m,k)}\cap K(\{k,m\},X)\subset T$ and hence
$K(\{k,m\},X)=T$, which proves the statement \eqref{pr302.2}.

To prove \eqref{pr302.3} we only have to observe that the (non-empty) set 
$X$ can  be chosen in  $2^{(l-1)^{n-l}}-1$ different ways.
\end{proof}

Our last construction corresponds to the case, described in 
Proposition~\ref{p18}\eqref{p18.2}.
For $k\in\{1,\dots,l\}$, $M\subset \{a_1,\dots,a_{k-1},a_{k+1},\dots,a_l\}$,
$|M|>1$, and $\varnothing\neq Y\subset\mathcal{V}^{(k)}$, let $L(k,M,Y)$ denote 
the union of all $\sqrt{\varepsilon}$, where $\varepsilon\in 
\mathcal{T}_n^{(l-1,k)}$ is an idempotent such that $\mathfrak{b}(\varepsilon)\in M$,
$\im(\varepsilon)\in Y$.

\begin{proposition}\label{pr303}
Let $k\in\{1,\dots,l\}$.
\begin{enumerate}[(a)]
\item\label{pr303.1}
For every $M$ and $Y$ as above the set $L(k,M,Y)$ is an  isolated subsemigroup 
of $\mathcal{T}_n^{(l-1,k)}$.
\item\label{pr303.2}
Every isolated subsemigroup of $\mathcal{T}_n^{(l-1,k)}$, which is not contained in
any $\mathcal{T}_n^{(l-1,k,m)}$, is equal to the semigroup 
$L(k,M,Y)$ for appropriate $M$ and $Y$.
\item\label{pr303.3}
The number of isolated subsemigroups of $\mathcal{T}_n^{(l-1,k)}$, which are 
not contained in any $\mathcal{T}_n^{(l-1,k,m)}$, equals
$(2^{l-1}-l)(2^{p}-1)$, where $p=\prod_{i\neq k}|A_i|$.
\end{enumerate}
\end{proposition}

\begin{proof}
One proves \eqref{pr303.1} using Lemma~\ref{l20}, Corollary~\ref{c202} and
Proposition~\ref{p18}\eqref{p18.2}. The statement
\eqref{pr303.3} follows from \eqref{pr303.2}
by a standard combinatorial calculation. We prove \eqref{pr303.2}.

Let  $T$ be an isolated subsemigroup of $\mathcal{T}_n^{(l-1,k)}$, which is 
not contained in  any $\mathcal{T}_n^{(l-1,k,m)}$. Set
\begin{displaymath}
M=\{\mathfrak{b}(\beta):\beta\in T\},\quad\quad\text{and}\quad\quad
Y=\{\im^{\mathrm{st}}(\beta):\beta\in T\}.
\end{displaymath}
Then $|M|>1$ by our assumption and thus $T\subset L(k,M,Y)$. The non-trivial
part is to prove that $L(k,M,Y)\subset T$. Lemma~\ref{l20}, Corollary~\ref{c202} and
Proposition~\ref{p18}\eqref{p18.2} tell us that for 
$\beta,\gamma\in T$ we have $\mathfrak{b}(\beta*\gamma)=\mathfrak{b}(\beta)$
and $\im^{\mathrm{st}}(\beta*\gamma)=\im^{\mathrm{st}}(\gamma)$. In particular,
this implies that for fixed $a_m\in M$ and $I\in Y$ there always exists some
idempotent, $\varepsilon\in T$, given by \eqref{idempnew}, such that 
$\mathfrak{b}(\varepsilon)=a_m$ and $\im(\varepsilon)=I$. 

As $|M|>1$, for the same reason  we can also fix another idempotent, 
$\tau\in T$, such that  $\mathfrak{b}(\tau)\neq a_m$ and $\im(\tau)=I$. Set 
\begin{displaymath}
\Lambda=\Lambda_{\varepsilon}=
E_1\cup\dots\cup E_{k-1}\cup E_{k+1}\cup\dots\cup E_l
\end{displaymath}
and recall that $a_i\in E_i$ for all $i=1,\dots,k-1,k+1,\dots,l$,
$a_k\in E_m$. We will need the following statement:

\begin{lemma}\label{sos}
Let $x\in N\setminus A$ and assume that $x\in E_s$. 
Let $t\in \{1,\dots,k-1,k+1,\dots,l\}$, $t\neq s$. 
Define 
\begin{displaymath}
\Lambda'=E'_1\cup\dots\cup E'_{k-1}\cup E'_{k+1}\cup\dots\cup E'_l
\end{displaymath}
in the following way: for $i=1,\dots,k-1,k+1,\dots,l$ set
\begin{displaymath}
E'_i=
\begin{cases}
E_i, & i\neq s,t;\\
E_s\setminus\{x\}, & i=s;\\     
E_t\cup\{x\}, & i=t.
\end{cases}
\end{displaymath}
Then $T$ contains some $\gamma$ such that
$\mathfrak{b}(\gamma)=a_m$, $\im^{\mathrm{st}}(\gamma)=I$ and
$\Lambda^{\mathrm{st}}_{\gamma}=\Lambda'$.
\end{lemma}

\begin{proof}
Let $\sigma$ be any permutation of $1,\dots,k-1,k+1,\dots,l$ such that 
$\sigma(t)=\mathfrak{nb}(\tau)$ and $\sigma(s)=\mathfrak{nb}(\varepsilon)$
(such $\sigma$ exists since $\mathfrak{b}(\varepsilon)\neq \mathfrak{b}(\tau)$).
Consider the element $\beta\in\mathcal{T}_n$ defined as follows:
\begin{gather*}
\beta(x)=a_k,\\
\beta(E_s\setminus\{x\})=e_{\sigma(s)},\\
\quad \beta(E_i)=e_{\sigma(i)}\quad \text{ for all }
i\in \{1,\dots,l\}\setminus\{k,s\}.
\end{gather*}
One checks that $\beta\in \sqrt{\varepsilon}\subset T$, thus
$\gamma=\beta*\tau\in T$. However, Corollary~\ref{c201} implies that
$\Lambda_{\gamma}=\Lambda'$. That $\mathfrak{b}(\gamma)=a_m$, 
$\im(\gamma)=I$ follows from Lemma~\ref{l20} and Corollary~\ref{c202}.
This completes the proof.
\end{proof}

Using Lemma~\ref{sos} inductively (starting from $\varepsilon$), one
shows that for every $\tilde{\Lambda}\in \mathcal{U}^{(k,m)}$ the semigroup
$T$ contains an idempotent, $\tilde{\varepsilon}$, such that 
$\Lambda_{\tilde{\varepsilon}}=\tilde{\Lambda}$, 
$\mathfrak{b}(\tilde{\varepsilon})=a_m$, 
and $\im(\tilde{\varepsilon})=I$. The statement of the proposition follows.
\end{proof}

Now we can gather in the harvest of the hard work above.

\begin{corollary}\label{cp22}
Let $T$ be an isolated subsemigroup of $(\mathcal{T}_n^{(l-1)},*)$. Then
either $T=H(k,m,X,Y)$ for some appropriate $k$, $m$, $X$, and $Y$, or
$T=K(\{k,m\},X)$ for some appropriate $\{k,m\}$ and $X$, 
or $T=L(k,M,Y)$ for some appropriate $k$, $M$, $Y$.
\end{corollary}

\begin{proof}
Let $T$ be an isolated subsemigroup of $(\mathcal{T}_n^{(l-1)},*)$.
Proposition~\ref{p18}\eqref{p18.4} implies that 
either $T\subset \mathcal{T}_n^{(l-1,k,m)}\cup\mathcal{T}_n^{(l-1,m,k)}$
for some appropriate $m$ and $k$ or
$T\subset \mathcal{T}_n^{(l-1,k)}$ for some appropriate $k$.
 
If $T\subset \mathcal{T}_n^{(l-1,k,m)}\cup\mathcal{T}_n^{(l-1,m,k)}$
and $T\subset \mathcal{T}_n^{(l-1,k)}$ at the same time, then
$T\subset \mathcal{T}_n^{(l-1,k,m)}$ and hence
$T=H(k,m,X,Y)$ for some appropriate $X$ and $Y$
by Proposition~\ref{pr301}.

If $T\subset \mathcal{T}_n^{(l-1,k,m)}\cup\mathcal{T}_n^{(l-1,m,k)}$
but $T\not\subset \mathcal{T}_n^{(l-1,k)}$ and
$T\not\subset \mathcal{T}_n^{(l-1,m)}$, then
$T=K(\{k,m\},X)$ for some appropriate $X$ by Proposition~\ref{pr302}.

Finally, if $T\subset \mathcal{T}_n^{(l-1,k)}$ but
$T\not \subset \mathcal{T}_n^{(l-1,k,m)}$ for any $m$, then
$T=L(k,M,Y)$ for some appropriate $M$ and $Y$ by Proposition~\ref{pr303}.
This completes the proof.
\end{proof}

\subsection{The case of stable rank at most $l-2$}\label{s4.3}

During this subsection we assume that $l>2$.

\begin{proposition}\label{pr25}
Let $T$ be an isolated subsemigroup of $(\mathcal{T}_n,*)$ such that
$T\not\subset \mathcal{T}_n^{(l)}\cup \mathcal{T}_n^{(l-1)}$. Then
$T\supset \mathcal{T}_n\setminus \mathcal{T}_n^{(l)}$.
\end{proposition}

To prove this we will need the following lemmas, in which without
loss of generality we assume that $e_i\in A_i$ for all $i=1,\dots,k$:

\begin{lemma}\label{l26}
Let $T$ be an isolated subsemigroup of $(\mathcal{T}_n,*)$.
Assume that $T$ contains an idempotent, $\varepsilon$, 
of rank $k$, given by \eqref{idemp}, where $1<k<l-1$. Then
$T$ contains an idempotent of rank $k-1$ or less.
\end{lemma}

\begin{proof}
Assume first that at least one of $E_i$'s contains at least three 
elements from $A$. Without loss of generality we may assume that
$E_1\supset\{a_1,a_l,a_{l-1}\}$. Consider the following elements:
\begin{gather*}
\beta=\left(
\begin{array}{ccccc}
E_1\setminus\{a_{l-1}\} & E_2 & \dots & E_k & a_{l-1} \\
e_1 & e_2 & \dots & e_k & a_l
\end{array}
\right),\\
\gamma=\left(
\begin{array}{ccccc}
E_1\setminus\{a_l\} & E_2 & \dots & E_k & a_l \\
e_1 & e_2 & \dots & e_k & a_{l-1}
\end{array}
\right),\\
\nu=\left(
\begin{array}{cccccc}
E_1\setminus\{a_l\} & E_2 & E_3 & \dots & E_k & a_l \\
e_1 &  a_l & e_3 & \dots & e_k & e_2
\end{array}
\right).
\end{gather*}
We have $\beta*\beta=\gamma*\gamma=\varepsilon$ and hence $\beta,\gamma\in T$.
Thus $T$ contains the element $\gamma*\beta$.
But $\nu*\nu=\gamma*\beta$, implying $\nu\in T$. At the same time
$\rank(\nu*\beta*\nu*\beta)=k-1$ and the statement follows.

Let us now assume that none of the idempotents in $T$ of rank at most
$l-2$ contains a block in the kernel having three or more elements
from $A$. 
Then, without loss of generality, we can assume that
$\{a_1,a_l\}\subset E_1$ and $\{a_2,a_{l-1}\}\subset E_2$.
Consider the following elements:
\begin{gather*}
\beta=\left(
\begin{array}{cccccc}
E_1\setminus\{a_{l}\} & E_2 & E_3 & \dots & E_k & a_{l} \\
e_2 & e_1 & e_3 & \dots & e_k & a_{l-1}
\end{array}
\right),\\
\gamma=\left(
\begin{array}{cccccc}
E_1 & E_2\setminus\{a_{l-1}\} & E_3 & \dots & E_k & a_{l-1} \\
e_2 & e_1 & e_3 &  \dots & e_k & a_{l}
\end{array}
\right),\\
\nu=\left(
\begin{array}{cccccc}
E_1 & E_2\setminus\{a_{l-1}\} & E_3 & \dots & E_k & a_{l-1} \\
e_2 &  e_1 & e_3 & \dots & e_k & a_{l-1}
\end{array}
\right).
\end{gather*}
We have $\beta*\beta=\gamma*\gamma=\varepsilon$ and hence
$\beta,\gamma\in T$, in particular, $\gamma*\beta\in T$. At the
same time $\nu*\nu=\gamma*\beta$, forcing $\nu\in T$.
But then $a_1$, $a_l$ and $a_{l-1}$ belong to the same block of
$\Lambda_{\nu*\varepsilon}$. Thus the element 
$\nu*\varepsilon*\nu*\varepsilon$, which is an idempotent in $T$ of 
rank at most $l-2$, has a block with  at least three elements from $A$. 
The obtained contradiction completes the proof.
\end{proof}

\begin{lemma}\label{l27}
Let $T$ be an isolated subsemigroup of $(\mathcal{T}_n,*)$.
Assume that $T$ contains an idempotent of rank $1$ and $l>2$. 
Then $T$ contains all idempotents
of $(\mathcal{T}_n,*)$, whose rank does not exceed $l-1$.
\end{lemma}

\begin{proof}
First let us show that $T$ contains all idempotents of
rank $1$. Without loss of generality we can assume that 
$T$ originally contains $\theta_1$ (the element with image
$\{1\}$) and  $1\in A_1$. Let $y\in A_3$. 
Consider the following elements:
\begin{gather*}
\beta=\left(
\begin{array}{cccc}
A_1 & A_2 &  \dots & A_l \\
1 & a_1 & \dots & a_{l-1}
\end{array}
\right),\\
\gamma=\left(
\begin{array}{cccccc}
A_1 & A_2 & A_3 & \dots & A_{l-1} & A_l \\
1 & y & a_4 &  \dots & a_l & 1
\end{array}
\right),\\
\nu=\left(
\begin{array}{cccccc}
A_1 & A_2 & A_3 & A_4 &  \dots  & A_l \\
y &  y & 1 & a_4 & \dots  & a_{l}
\end{array}
\right).
\end{gather*}
We have $\beta^{*n}=\gamma^{*n}=\theta_1$ and hence
$\beta,\gamma\in T$. In particular, $\beta*\gamma\in T$.
But $\nu*\nu=\beta*\gamma$, hence $\nu\in T$. At the same time
$\theta_1*\nu=\theta_y$ and thus $\theta_y\in T$. This actually shows that 
$\theta_y\in T$ for all $y\not\in A_1$. But then, taking any idempotent of 
rank one with the image outside $A_1$, from the above it also follows that 
$\theta_y\in T$ for all $y\in A_1$. Observe that the construction above 
does not work in the cases $l=1,2$.

Now let $\varepsilon$ be an idempotent of $(\mathcal{T}_n,*)$
of rank $k\leq l-1$ given by \eqref{idemp}. Without loss of generality 
we may assume $a_l\in E_1$. Consider the following elements:
\begin{gather*}
\beta=\left(
\begin{array}{cccccc}
E_1 & E_2 & E_3 & \dots & E_{k-1} & E_k \\
e_1 & a_3 & a_4 & \dots & a_k & a_{l}
\end{array}
\right),\\
\gamma=\left(
\begin{array}{cccccc}
N\setminus\{a_3,\dots,a_k,a_l\} & a_3 & a_4 & \dots & a_k & a_{l} \\
e_1 & e_2 & e_3 &  \dots & e_{k-1} & e_k
\end{array}
\right).
\end{gather*}
From the first part of the proof we have $\theta_{e_1}\in T$.
We further have $\beta^{*n}=\gamma^{*n}=\theta_{e_1}$ and hence
$\beta,\gamma\in T$. At the same time
$\beta*\gamma=\varepsilon$ also must belong to $T$. This completes the
proof.
\end{proof}

Now we are ready to prove Proposition~\ref{pr25}:

\begin{proof}[Proof of Proposition~\ref{pr25}]
If $T\not\subset \mathcal{T}_n^{(l)}\cup \mathcal{T}_n^{(l-1)}$,
it must contain an idempotent of rank $l-2$ or less. Applying
Lemma~\ref{l26} inductively we obtain that $T$ must contain
an idempotent of rank $1$. From Lemma~\ref{l27} it now follows that
$T$ contains all idempotents of rank $\leq l-1$.
Since $T$ is isolated, we derive that 
$T\supset \mathcal{T}_n\setminus \mathcal{T}_n^{(l)}$ and complete the proof.
\end{proof}

\subsection{General classification of isolated subsemigroups}\label{s4.4}

After the hard work in Subsections~\ref{s4.1}--\ref{s4.3} we can now 
relatively easily complete the classification of all isolated 
subsemigroups of  $(\mathcal{T}_n,*)$. We write
$\mathcal{T}_n$ as the disjoint union of the sets
$\mathcal{T}_n^{(l)}$, $\mathcal{T}_n^{(l-1)}$, and
$Z=\mathcal{T}_n\setminus(\mathcal{T}_n^{(l)}\cup \mathcal{T}_n^{(l-1)})$,
and for every collection of these sets we will classify all isolated 
subsemigroups, which have non-trivial intersections exactly with the 
elements of this collection. 

\begin{theorem}\label{tmain1}
Let $T$ be an isolated subsemigroup of $(\mathcal{T}_n,*)$. Then $T$ is
one of the following:
\begin{enumerate}[(i)]
\item\label{tmain1.1} $F(X,Y)$ for appropriate $X\subset \mathcal{X}$ and 
$Y\subset \mathcal{Y}$
(see Subsection~\ref{s4.1});
\item\label{tmain1.20} $H(k,m,X,Y)$ for appropriate  $k,m\in\{1,\dots,m\}$,
$k\neq m$, $X\subset \mathcal{U}^{(k,m)}$, and $Y\subset \mathcal{V}^{(k)}$, if $l>1$
(see Subsection~\ref{s4.2});
\item\label{tmain1.21} $K(\{k,m\},X)$ for appropriate  
$\{k,m\}\subset\{1,\dots,l\}$, $k\neq m$, and  
$X\subset \mathcal{U}^{(k,m)}$, if $l>1$ (see Subsection~\ref{s4.2});
\item\label{tmain1.22}
$L(k,M,Y)$ for appropriate $k\in\{1,\dots,l\}$, $Y\subset \mathcal{V}^{(k)}$, and
\begin{displaymath}
M\subset\{a_1,\dots,a_{k-1},a_{k+1},\dots,a_l\},\quad\quad |M|>1,
\end{displaymath}
if $l>1$ (see Subsection~\ref{s4.2});
\item\label{tmain1.3} $\mathcal{T}_n\setminus \mathcal{T}_n^{(l)}$, if $l>2$;
\item\label{tmain1.5} $F(X,Y)\cup (\mathcal{T}_n\setminus \mathcal{T}_n^{(l)})$,
for appropriate  $X\subset \mathcal{X}$ and $Y\subset \mathcal{Y}$, if $l>1$;
\end{enumerate}
All semigroups in the above list are different, in particular, the list above
gives a complete classification of isolated subsemigroup of $(\mathcal{T}_n,*)$.
\end{theorem}

\begin{proof}
We have to consider the following cases:

\noindent
{\bf Case 1:} $T$ intersects only with $\mathcal{T}_n^{(l)}$. In this case
Corollary~\ref{ccc123} gives \eqref{tmain1.1}.

\noindent
{\bf Case 2:} $T$ intersects only with $\mathcal{T}_n^{(l-1)}$. This makes
sense only if $l>1$. In this case
Corollary~\ref{cp22} gives \eqref{tmain1.20}, 
\eqref{tmain1.21} and \eqref{tmain1.22}.

\noindent
{\bf Case 3:} $T$ intersects only with $Z$. This makes
sense only if $l>2$, when it is not possible by
Proposition~\ref{pr25}.

\noindent
{\bf Case 4:} $T$ intersects only with $Z$ and $\mathcal{T}_n^{(l)}$. 
This makes sense only if $l>2$, when it is not possible by 
Proposition~\ref{pr25}.

\noindent
{\bf Case 5:} $T$ intersects only with $Z$ and $\mathcal{T}_n^{(l-1)}$. 
This makes sense only for $l>2$.
Since $\mathcal{T}_n\setminus \mathcal{T}_n^{(l)}$ is an ideal and the
complement of an isolated subsemigroup, it is an isolated subsemigroup
itself. Hence in this case Proposition~\ref{pr25} gives \eqref{tmain1.3}.

\noindent
{\bf Case 6:} $T$ intersects only with $\mathcal{T}_n^{(l)}$ and 
$\mathcal{T}_n^{(l-1)}$.  This makes sense only for $l>1$.
Let $\varepsilon\in T$ be an idempotent of 
rank $l$ and $\nu\in T$ be an idempotent of rank
$l-1$. Corollary~\ref{c7} implies that for any permutation on $A$ there exists 
an element from $G(\varepsilon)$ which induces on $A$ this fixed permutation.
In particular, if $l>2$ then, multiplying $\nu$ from the left and from the right 
with the elements from $G(\varepsilon)$, we can change $\mathfrak{b}(\nu)$ and
$\mathfrak{t}(\nu)$ in an arbitrary way. Thus we can find some
$\nu'\in T$ such that neither 
$\mathfrak{d}(\nu)=\mathfrak{d}(\nu')$ nor
$\mathfrak{t}(\nu)=\mathfrak{t}(\nu')$ holds. From
Proposition~\ref{p18}\eqref{p18.4} it follows that such $T$ must intersect with
$Z$, a contradiction. Hence for $l>2$ this case is not possible. 

Let us now assume that $l=2$. Then $T\cap \mathcal{T}_n^{(2)}=F(X,Y)$
for appropriate $X$ and $Y$ because of Corollary~\ref{ccc123}. Further,
Proposition~\ref{pr302}\eqref{pr302.2} implies that
$T\cap \mathcal{T}_n^{(1)}=\mathcal{T}_n^{(1)}=
\mathcal{T}_n\setminus \mathcal{T}_n^{(2)}$, which
gives the case \eqref{tmain1.5} for $l=2$. 

\noindent
{\bf Case 7:} $T$ intersects with $\mathcal{T}_n^{(l)}$,  
$\mathcal{T}_n^{(l-1)}$ and $Z$. This makes sense only if $l>2$.
In this case Proposition~\ref{pr25} implies
that $T$ contains the ideal $\mathcal{T}_n\setminus \mathcal{T}_n^{(l)}$. The 
intersection with $\mathcal{T}_n^{(l)}$ should be one of $F(X,Y)$ and since
$\mathcal{T}_n\setminus \mathcal{T}_n^{(l)}$ is an ideal, it follows that
the union $F(X,Y)\cup (\mathcal{T}_n\setminus \mathcal{T}_n^{(l)})$ is a semigroup
for all possible $X$ and $Y$. It is isolated by Proposition~\ref{pr3}. This gives
\eqref{tmain1.5}.

Since all possible cases are considered, the statement follows.
\end{proof}

\begin{corollary}\label{skvert}
Let $\varepsilon\in E(\mathcal{T}_n,*)$. Then $\sqrt{\varepsilon}$ is a
subsemigroup if and only if $\rank(\varepsilon)=l$ or $\rank(\varepsilon)=l-1$.
\end{corollary}

\begin{corollary}\label{number}
Set $p=\prod_{i=1}^l|A_i|$ and for $k=1,\dots,l$ set $p_k=\prod_{i\neq k}|A_i|$. Then
\begin{enumerate}[(a)]
\item\label{number.2} for $l=1$ the semigroup $(\mathcal{T}_n,*)$ contains 
$2^n-1$ isolated subsemigroups;
\item\label{number.3} for $l=2$ the semigroup $(\mathcal{T}_n,*)$ contains 
$2^{2^{n-2}+p}-2^{2^{n-2}}-2^{p}+2^n$ isolated subsemigroups;
\item\label{number.1}
for $l>2$ the semigroup $(\mathcal{T}_n,*)$ contains
\begin{multline*}
\sum_{k=2}^l\sum_{m=1}^{k-1}(2^{(l-1)^{n-l}}-1)(2^{p_k}-1)(2^{p_m}-1)+\\
+\sum_{k=1}^l l\cdot (2^{(l-1)^{n-l}+p_k}-2^{(l-1)^{n-l}}-2^{p_k}+1)
+\sum_{k=1}^l(2^{l-1}-l)(2^{p_k}-1)+\\
+ 2\cdot (2^{l^{n-l}+p}-2^{l^{n-l}}-2^{p}+1)+1
\end{multline*}
isolated subsemigroups.
\end{enumerate}
\end{corollary}

\begin{proof}
All statements 
follow from Theorem~\ref{tmain1}, Proposition~\ref{pr14}\eqref{pr14.4},
Proposition~\ref{pr301}\eqref{pr301.5}, Proposition~\ref{pr302}\eqref{pr302.3},
and Proposition~\ref{pr303}\eqref{pr303.3} using standard combinatorial computations.
\end{proof}

\section{Classification of other types of isolated subsemigroups for $l>2$}\label{s5}

\begin{theorem}\label{tmain2}
Let $l>2$. Then
\begin{enumerate}
\item\label{tmain2.1}
\begin{multline*}
\mathfrak{CI}(\mathcal{T}_n,*)=\{
\mathcal{T}_n,\mathcal{T}_n^{(l)},\mathcal{T}_n\setminus \mathcal{T}_n^{(l)}
\}\cup\\
\cup\{ F(X,\mathcal{Y}),
F(X,\mathcal{Y})\cup(\mathcal{T}_n\setminus \mathcal{T}_n^{(l)}):
X\subset \mathcal{X}, X\neq \varnothing, \mathcal{X}\}\cup\\
\cup\{F(\mathcal{X},Y),
F(\mathcal{X},Y)\cup(\mathcal{T}_n\setminus \mathcal{T}_n^{(l)}):
Y\subset \mathcal{Y}, Y\neq \varnothing, \mathcal{Y}
\}.
\end{multline*}
\item\label{tmain2.2}
$\mathfrak{RC}(\mathcal{T}_n,*)=\{\mathcal{T}_n\}\cup 
\{F(\mathcal{X},Y): Y\subset \mathcal{X}, Y\neq \varnothing\}$.
\item\label{tmain2.3}
$\mathfrak{LC}(\mathcal{T}_n,*)=\{\mathcal{T}_n\}\cup \{ F(X,\mathcal{Y}):
X\subset \mathcal{X}, X\neq \varnothing\}$.
\item\label{tmain2.4}
$\mathfrak{C}(\mathcal{T}_n,*)=\{\mathcal{T}_n,\mathcal{T}_n^{(l)}\}$.
\end{enumerate}
\end{theorem}

\begin{proof}
We prove \eqref{tmain2.1} and leave the proof of the other statements 
to the reader. Because of \eqref{br9}, to prove \eqref{tmain2.1} it is enough
to check which of the semigroups listed in Theorem~\ref{tmain1} are completely
isolated, that is have a complement, which is a semigroup itself.

$\overline{F(X,Y)}$ is a subsemigroup if and only if
$\mathcal{T}_n^{(l)}\setminus \overline{F(X,Y)}$ is a subsemigroup. Using
Proposition~\ref{pr14}\eqref{pr14.3} we easily obtain that the latter is equivalent
to $X=\mathcal{X}$ or $Y=\mathcal{Y}$, that is $F(X,Y)\in \mathfrak{CI}(\mathcal{T}_n,*)$
if and only if $X=\mathcal{X}$ or $Y=\mathcal{Y}$.

If $l>2$ and $H(k,m,X,Y)$ is completely isolated, then the complement
$\overline{H(k,m,X,Y)}$ is completely isolated as well, in particular,
it is isolated. However, it contains $\mathcal{T}_n^{(l)}$ and intersects with
$\mathcal{T}_n^{(l-1)}$, which, by Theorem~\ref{tmain1}, implies that 
$\overline{H(k,m,X,Y)}=\mathcal{T}_n$, a contradiction. Hence
$H(k,m,X,Y)$ is not completely isolated. Same arguments also imply that
neither the subsemigroup $K(\{k,m\},X)$ nor $L(k,M,Y)$ is  
completely isolated.

$\mathcal{T}_n\setminus \mathcal{T}_n^{(l)}$ is completely isolated as the
complement to the completely isolated subsemigroup $\mathcal{T}_n^{(l)}=
F(\mathcal{X},\mathcal{Y})$.

Finally, the complement to 
$F(X,Y)\cup(\mathcal{T}_n\setminus \mathcal{T}_n^{(l)})$ belongs to
$\mathcal{T}_n^{(l)}$ and is completely isolated if and only if
$F(X,Y)\cup(\mathcal{T}_n\setminus \mathcal{T}_n^{(l)})$ is. The above classification of
completely isolated subsemigroups of $\mathcal{T}_n^{(l)}$ implies that
$F(X,Y)\cup(\mathcal{T}_n\setminus \mathcal{T}_n^{(l)})$ is completely isolated
if and only if $X=\mathcal{X}$ or $Y=\mathcal{Y}$. This completes the proof.
\end{proof}

\section{Exceptional cases $l=1,2$ and an application to the case $l=n$}\label{s6}

In this section we separately formulate the main statement for the
extreme case $l=n$, that is for the case of the classical $\mathcal{T}_n$,
since this is clearly a case of independent interest. We also separately consider
the exceptional case $l=1,2$, and in each of these cases give a short formulation
of Theorem~\ref{tmain1} and present an analogue of Theorem~\ref{tmain2}.

\subsection{The case of the classical $\mathcal{T}_n$}\label{s6.1}

For $m,k\in N$, $m\neq k$, let $\varepsilon_{m,k}$ be defined as
in \eqref{p18.eq.1}. In this subsection we assume that
$\alpha=\mathrm{id}$, that is $(\mathcal{T}_n,*)=\mathcal{T}_n$.
Recall (see Corollary~\ref{c7}) that the maximal subgroup
$G(\varepsilon_{m,k})$ consists of all $\beta\in \mathcal{T}_n$ such that
$\beta(m)=\beta(k)$ and which induce a permutation on $N\setminus\{k\}$.

\begin{proposition}\label{clast1}
\begin{enumerate}[(a)]
\item\label{clast1.1}  $\mathfrak{I}(\mathcal{T}_1)=
\mathfrak{CI}(\mathcal{T}_1)=\mathfrak{LC}(\mathcal{T}_1)=
\mathfrak{RC}(\mathcal{T}_1)=\mathfrak{C}(\mathcal{T}_1)= \{\mathcal{T}_1\}$.
\item\label{clast1.2} 
\begin{enumerate}[(i)]
\item\label{clast1.2.1}
$\mathfrak{I}(\mathcal{T}_2)=
\{\mathcal{T}_2,\mathcal{S}(\{1,2\}),\mathcal{T}_2\setminus\mathcal{S}(\{1,2\}),
G(\varepsilon_{1,2}),G(\varepsilon_{2,1}) \}$.
\item\label{clast1.2.2}
$\mathfrak{CI}(\mathcal{T}_2)=
\left\{ \mathcal{T}_2,\mathcal{S}(\{1,2\}), 
\mathcal{T}_2\setminus\mathcal{S}(\{1,2\})
\right\}$.
\item\label{clast1.2.3}
$\mathfrak{LC}(\mathcal{T}_2)=\mathfrak{RC}(\mathcal{T}_2)=
\mathfrak{C}(\mathcal{T}_2)= \{\mathcal{S}(\{1,2\}),\mathcal{T}_2\}$
\end{enumerate}
\item\label{clast1.3} Let $n>2$. Then
\begin{enumerate}[(i)]
\item\label{clast1.3.1}
\begin{multline*}
\mathfrak{I}(\mathcal{T}_n)=
\{\mathcal{T}_n,\mathcal{S}(N),\mathcal{T}_n\setminus\mathcal{S}(N)\}
\cup\\\cup
\left\{\bigcup_{m\in M}G(\varepsilon_{m,k}): k\in N, \varnothing\neq 
M\subset N\setminus\{k\}
\right\}
\cup\\\cup
\left\{G(\varepsilon_{m,k})\cup 
G(\varepsilon_{k,m}): m,k\in N, m\neq k \right\}.
\end{multline*}
\item\label{clast1.3.2}
$\mathfrak{CI}(\mathcal{T}_n)=
\left\{ \mathcal{T}_n,\mathcal{S}(N), \mathcal{T}_n\setminus\mathcal{S}(N)
\right\}$.
\item\label{clast1.3.3}
$\mathfrak{LC}(\mathcal{T}_n)=\mathfrak{RC}(\mathcal{T}_n)=
\mathfrak{C}(\mathcal{T}_n)= \{\mathcal{S}(N),\mathcal{T}_n\}$
\end{enumerate}
\end{enumerate}
\end{proposition}

\begin{proof}
The statement \eqref{clast1.3} follows from Theorem~\ref{tmain1} and 
Theorem~\ref{tmain2}. The statements \eqref{clast1.1} and \eqref{clast1.2}
are obvious.
\end{proof}

\subsection{The case when the sandwich element has rank $1$}\label{s6.2}

\begin{proposition}\label{clast2}
Assume that $\alpha=\theta_1$. Then 
\begin{enumerate}[(a)]
\item\label{clast2.1} $E(\mathcal{T}_n,*)=\{\theta_i:i=1,\dots,n\}$.
\item\label{clast2.2} 
\begin{displaymath}
\mathfrak{I}(\mathcal{T}_n,*)=
\mathfrak{CI}(\mathcal{T}_n,*)=\mathfrak{RC}(\mathcal{T}_n,*)=
\left\{ \bigcup_{i\in X}\sqrt{\theta_i}: \varnothing\neq X\subset N\right\}.
\end{displaymath}
\item\label{clast2.3} $\mathfrak{LC}(\mathcal{T}_n,*)=\mathfrak{C}(\mathcal{T}_n,*)=
\{\mathcal{T}_n\}$. 
\end{enumerate}
\end{proposition}

\begin{proof}
The statement \eqref{clast2.1} follows from Theorem~\ref{t6}\eqref{t6.1}.
The statement about $\mathfrak{I}(\mathcal{T}_n,*)$ follows from
Theorem~\ref{tmain1}. Lemma~\ref{l11}\eqref{l11.2} implies  that
in the case $\alpha=\theta_1$ every isolated subsemigroup of 
$\mathfrak{I}(\mathcal{T}_n,*)$ is right convex, in particular, 
is completely isolated. This proves \eqref{clast2.2}. 
Again Lemma~\ref{l11}\eqref{l11.2} implies  that
in the case $\alpha=\theta_1$  every proper subsemigroup of
$\mathfrak{I}(\mathcal{T}_n,*)$ can neither be left convex nor convex.
This proves \eqref{clast2.3} and completes the proof.
\end{proof}

\subsection{The case when the sandwich element has rank $2$}\label{s6.3}

\begin{proposition}\label{clast3}
Assume that $l=2$. Then 
\begin{enumerate}[(a)]
\item\label{clast3.1}
\begin{multline*}
\mathfrak{I}(\mathcal{T}_n,*)=\{F(X,Y),F(X,Y)\cup \mathcal{T}_n^{(1)}
: X\subset \mathcal{X},
Y\subset \mathcal{Y};X,Y\neq \varnothing\}\cup\\
\cup \left\{\bigcup_{i\in X}\sqrt{\theta_i}:
X\subset N,X\neq \varnothing\right\}.
\end{multline*} 
\item\label{clast3.2} 
\begin{multline*}
\mathfrak{CI}(\mathcal{T}_n,*)=\left\{\mathcal{T}_n,
\mathcal{T}_n^{(2)},\mathcal{T}_n^{(1)}\right\}\cup\\ \cup
\left\{F(\mathcal{X},Y),F(\mathcal{X},Y)\cup 
\mathcal{T}_n^{(1)}:Y\subset \mathcal{Y}, Y\neq \varnothing,\mathcal{Y}
\right\}\cup\\
\cup \{F(X,\mathcal{Y}),F(X,\mathcal{Y})\cup \mathcal{T}_n^{(1)}:
X\subset \mathcal{X}, X\neq\varnothing, \mathcal{X}\}.
\end{multline*} 
\item\label{clast3.3} 
$\mathfrak{RC}(\mathcal{T}_n,*)=\left\{\mathcal{T}_n\right\}
\cup \{F(X,\mathcal{Y}): X\subset \mathcal{X}, X\neq\varnothing\}$.
\item\label{clast3.4} 
$\mathfrak{LC}(\mathcal{T}_n,*)=\left\{\mathcal{T}_n\right\}\cup 
\left\{F(\mathcal{X},Y):Y\subset \mathcal{Y}, Y\neq \varnothing\right\}$.
\item\label{clast3.5} 
$\mathfrak{C}(\mathcal{T}_n,*)=\{\mathcal{T}_n,\mathcal{T}_n^{(2)}\}$.
\end{enumerate}
\end{proposition}

\begin{proof}
The statement \eqref{clast3.1} follows from 
Theorem~\ref{tmain1}\eqref{tmain1.1}--\eqref{tmain1.22},\eqref{tmain1.5}.

To prove \eqref{clast3.2} we should go through the list from \eqref{clast3.1}
and see for which semigroups there the complement is again a semigroup, that is
again occurs in the list from \eqref{clast3.1}. This is a  straightforward
calculation.

To prove \eqref{clast3.3} and \eqref{clast3.4} we should go through the list 
from \eqref{clast3.2} and see for which semigroups there the complement
is a left or a right ideal respectively. This is again a  straightforward
calculation.

Finally, a semigroup is convex if and only if it is both left convex and
right convex. Hence the list from \eqref{clast3.5} is just the intersection
of the lists from \eqref{clast3.3} and \eqref{clast3.4}. This completes the
proof.
\end{proof}

\vspace{0.1cm}

\begin{center}
\bf Acknowledgments
\end{center}

An essential part of the  paper was written during the visit of the 
second author to Uppsala University, which was supported by 
The Swedish Institute. The financial support of The Swedish Institute 
and the hospitality of Uppsala University are gratefully 
acknowledged. For the first author the research was partially 
supported by The Swedish Research Council. We are grateful to the
referee for pointing out an error in the original version of the
paper and for many very useful suggestions and comments.

\vspace{0.1cm}

\noindent
V.M.: Department of Mathematics, Uppsala University, Box 480, 
SE 751 06, Uppsala, SWEDEN, e-mail: {\em mazor\symbol{64}math.uu.se},\\ web: 
``http://www.math.uu.se/$\tilde{\hspace{2mm}}$mazor''
\vspace{0.2cm}

\noindent
G.T.: Department of Mechanics and Mathematics, Kyiv Taras Shevchenko 
University, 64, Volodymyrska st., 01033, Kyiv, UKRAINE,\\ e-mail: 
{\em gtsyaputa\symbol{64}univ.kiev.ua}

\end{document}